\documentclass[a4paper,12pt,twoside]{article}
\usepackage[latin1]{inputenc}
\usepackage[T1]{fontenc}
\usepackage{amssymb}
\usepackage{dochier}
\usepackage[final]{lpretex}
\usepackage{title}
\usepackage{a4}
\usepackage{amscd}
\usepackage{caption}
\input xy
\xyoption{all}

\EndOfDump

\def\DHpartformat#1{(#1) }
\def\DHrefpart#1{(\DHRefpart{#1})}
\LBsetstyleof{partno}{arabic}

\def\i {\item}

\let\define\def

\def\A {{\mathbb A}}  \def\C {{\mathbb C}}
  
\def\GG {{\mathbb G}} \def\H{{\mathbb H}}  
  
  \def\P {{\mathbb P}} 
\def\Q {{\mathbb Q}} 
 
  \def\X {{\mathbb X}}
\def\Z {{\mathbb Z}} 

\define \n {\mathbb N}
\define \q {\mathbb Q}
\define \PP {\mathbb P}

\def\sA {{\Cal A}} \def\sB {{\Cal B}} \def\sC {{\Cal C}}
 \def\sE {{\Cal E}} \def\sF {{\Cal F}}
\def\sG {{\Cal G}} \def\sH {{\Cal H}} \def\sI {{\Cal I}}
\def\sJ {{\Cal J}}  
\def\sM {{\Cal M}} \def\sN {{\Cal N}} \def\sO {{\Cal O}}
 
\def\sS {{\Cal S}} \def\sT {{\Cal T}} 
\def\sV {{\Cal V}}  \def\sX {{\Cal X}}
\def\sY {{\Cal Y}}

\def\sZ {{\Cal Z}}
\define \cN {\Cal N}
\define \cf {\Cal F}
\define \cg {\Cal G}
\define \cE {\Cal E}
\define \ce {\Cal E}
\define \cc {\Cal C}
\define \cV {\Cal V}
\define \cA {\Cal A}
\define \cK {\Cal K}
\define \cO {\Cal O}
\define \cF {\Cal F}
\define \cn {\Cal N}
\define \cI {\Cal I}
\define \sP {\Cal P}
\define \sEll {\Cal{Ell}}
\define \sJE {\Cal{JE}}
\define \sse{\Cal{SE}}

\define \sGJE {\Cal{GJE}}
\define \sHyp {\mathcal{HYP}}
\def\a {\alpha} \def\b {\beta} \def\g {\gamma}  
 \def\z{\zeta}
  
\def\s {\sigma}

\define \x {\xi}
\define \y {\eta}
\define \G {\Gamma}
\define \r {\rho}
\define \w {\omega}

\def\tX {\widetilde X}

\def \tC {\widetilde C}

\def \trho {\tilde {\rho}}

\def\tP {\widetilde P}
\define \tH {\widetilde H}
\define \tG {\widetilde{G}}
\define \tW {\widetilde W}
\define \tF {\widetilde F}
\define \tm {\tilde m}
\define \St {\widetilde S}
\define \Xt {\widetilde X}
\define \tS {\widetilde S}
\define \tpsi {\tilde \psi}
\define \tL {\widetilde L}
\define \tE {\widetilde E}
\define \tl {\tilde l}
\define \tA {\widetilde A}
\define \tom {\tilde\omega}
\define \tT {\widetilde T}
\define \tB {\widetilde B}
\define \tf {\tilde f}
\define \tsA {\widetilde{\sA}}

\define \tM {\widetilde M}
\define \tpsi {\widetilde{\psi}}
\define \trho {\widetilde{\rho}}
\define \tR {\widetilde R}
\define \tp {\widetilde p}
\define \tq {\widetilde q}
\define \tc {\tilde c}
\define \tsF {\widetilde {\sF}}
\define \tsM {\widetilde {\sM}}
\define \tii {\tilde i}
\define \tx {\tilde x}
\define \tg {\tilde g}
\define \tw {\tilde w}
\define \tz {\tilde z}
\define \ta {\widetilde\alpha}

\define \tsje{\widetilde{\sJE}}

\define \ts{\widetilde{\sigma}}

\define \bD {\overline{D}}
\define \bG {\overline{G}}
\define \bI {\overline{I}}
\define \bK {\overline{K}}
\define\bsM{\overline{\sM}}
\define\bsm{\bsM}
\define\sje{\sJE}
\define\bsell{\bsEll}
\define \bV {\overline{V}}
\define \bX {\overline{X}}
\define \bY {\overline{Y}}
\define \btau {\overline{\tau}}

\define\tDelta{\widetilde{\Delta}}
\define\bnabla{\overline{\nabla}}

\def\pd {\partial}

\def \Dx1 {\frac{\pd}{{\pd} x_1}}
\def \Dy1 {\frac{\pd}{{\pd} y_1}}
\def \Dz1 {\frac{\pd}{{\pd} z_1}}
\def \Dx2 {\frac{\pd}{{\pd} x_2}}
\def \Dy2 {\frac{\pd}{{\pd} y_2}}
\def \Dz2 {\frac{\pd}{{\pd} z_2}}

\def\q {\quad} 

\def\Mapdiagr#1{\nearrow\rlap{$\lower 5pt\vbox{{\hbox{$\mkern
-15mu\scriptstyle#1$}}}$}} 

\def\Mapdiagl#1{\llap{$\lower 5pt\vbox{{\hbox{$\scriptstyle#1\mkern
-15mu$}}}$}\searrow} 

\def\Mapswr#1{\swarrow\rlap{$\lower 5pt\vbox{{\hbox{$\mkern
-15mu\scriptstyle#1$}}}$}}              

\def\Mapnwl#1{\nwarrow\rlap{$\lower 5pt\vbox{{\hbox{$\mkern
-15mu\scriptstyle#1$}}}$}}

\def\i.e#1#2#3{\mathrel{\smash{\mathop{#2}\limits^{#1}_{#3}}}}

\def \inj {\hookrightarrow}

\define \Rhook {\hookrightarrow}

\def \half {\raise1pt\hbox{$\scriptstyle
        \frac{1}{2}\displaystyle$}}

\def \x{{\sl X}\llap{$\mkern -2mu {\scriptstyle -}$}}

\def \Symm {\operatorname{Sym}}
\def \Res {\operatorname{Res}}

\def \Pic {\operatorname{Pic}}

\define \Kod {\operatorname{Kod}}
\define \dimension {\operatorname{dim}}
\define \codim {\operatorname{codim}}
\define \contr {\operatorname{contr}}
\define \rk {\operatorname{rank}}
\define \Im {\operatorname {Im}}
\define \Mor {\operatorname{Mor}}
\define \Cl {\operatorname{Cl}}
\define \Hilb {\operatorname{Hilb}}
\define \degree {\operatorname{deg}}
\define \mult {\operatorname{mult}}
\define \Aut {\operatorname{Aut}}
\define \NS {\operatorname{NS}}
\define \Gal {\operatorname{Gal}}
\define \ch {\operatorname{char}}
\define \Jac {\operatorname{Jac}}
\define \Km {\operatorname{Km}}
\define \Sec {\operatorname{Sec}}
\define \Stab {\operatorname{Stab}}
\define \Br {\operatorname{Br}}
\define \Inv {\operatorname {Inv}}
\define \tr {\operatorname{tr}}
\define \Frob {\operatorname{Frob}}
\define \Symn {\operatorname{Symm}^n}
\define \Ev {\sE^\vee}
\define \ordp {\operatorname{ord}_p}
\define \Supp {\operatorname{Supp}}
\define \Ann {\operatorname{Ann}}
\define \disc {\operatorname{disc}}
\define \lie {\operatorname{lie}}
\define \embdim {\operatorname{embdim}}

\def\Def{\operatorname{Def}}

\def\Ram{\operatorname{Ram}}

\define\Ad{\operatorname{Ad}}

\def\bsEll{\overline{\sE\ell\ell}}

\define\nbd{neighbourhood }

\define\Tors{\operatorname{Tors}}
\define\ut{\underline{t}}
\define\u0{\underline{0}}
\define\uq{\underline{q}}
\define\uom{\underline{\omega}}

\def\hod#1#2#3#4{\ensuremath{ if#30 H^{#2}({#1},{\cal O}_{#1}) \else 
 H^{#2}(#1,\Omega^{#3}\if\relax{#4}\relax_{#1}\else _{#1/#4}\fi)\fi}}

\begin{document}
\title[Periods of elliptic surfaces]
{Generic Torelli and local Schottky theorems for  
 Jacobian elliptic surfaces}
\author{N. I. Shepherd-Barron}
\address{King's College,
Strand,
London WC2R 2LS,
U.K.}
\email{Nicholas.Shepherd-Barron@kcl.ac.uk}
\maketitle
\begin{abstract}
  Suppose that $f:X\to C$ is a general
  Jacobian elliptic surface
  over $\C$
  of irregularity $q$
  and positive geometric genus $h$.
  Assume that $10 h>12(q-1)$, that $h>0$
  and let $\bsell$ denote the stack of generalized elliptic curves.
  \begin{enumerate}
  \item The moduli stack $\sje$ of such surfaces
    is smooth at the point $X$ and its tangent space $T$ there
    is naturally a direct sum of lines
    $(v_a)_{a\in Z}$, where $Z\subset C$ is the
    ramification locus of the classifying morphism
    $\phi:C\to\bsell$ that corresponds to $X\to C$.
  \item For each $a\in Z$
    the map $\bnabla_{v_a}:H^{2,0}(X)\to H^{1,1}_{prim}(X)$
    defined by the derivative $per_*$ of the period map $per$
    is of rank one.
    Its image
    is a line $\C[\eta_a]$ and its kernel is
    $H^0(X,\Omega^2_X(-E_a))$,
    where $E_a=f^{-1}(a)$.
  \item The classes $[\eta_a]$ form an orthogonal
    basis of $H^{1,1}_{prim}(X)$ and $[\eta_a]$ is represented
    by a meromorphic $2$-form $\eta_a$ in $H^0(X,\Omega^2_X(2E_a))$
    of the second kind.
  \item We prove a local Schottky theorem; that is,
    we give a description of $per_*$ in terms of a certain
    additional structure on the vector bundles
    that are involved.

    Assume further that $8h>10(q-1)$ and that $h\ge q+3$.
    
  \item Given the period point $per(X)$ of $X$ that classifies the
    Hodge structure on the primitive cohomology
    $H^2_{prim}(X)$ and the image of $T$ under $per_*$
    we recover $Z$ as a subset of $\P^{h-1}$ and then, by
    quadratic interpolation, the curve $C$.
  \item We prove a generic Torelli theorem for
    these surfaces.
  \end{enumerate}
  Everything relies on the construction, via
  certain kinds of
  Schiffer variations\footnote{In an earlier version of this paper
    we used variations constructed by Fay \cite{F1}, \cite{F2}. However, Schiffer
    variations
    are slightly more powerful.} of curves,
    of certain 
    variations of $X$ for which $per_*$ can be calculated. 
  \end{abstract}
\begin{section}{Introduction}\label{intro}
  Suppose that $\sM$ is a moduli stack of smooth projective varieties over $\C$
  and and that $per:\sM\to\sP=D/\G$ is a corresponding period map. The
  derivative of $per$ is a homomorphism
  $$per_*:T_{\sM}\to per^*T_{\sP}.$$
  The \emph{local Torelli problem}
  is that of describing the kernel of this
  homomorphism and
  the \emph{local Schottky problem}
  is the problem of describing its image.
   We say that the
  \emph{local Torelli theorem}
  holds at a point $x$ of $\sM$ if the derivative $per_*$ of $per$
  is injective at $x$ and that the \emph{generic local Torelli theorem}
  holds if it holds at every generic point of $\sM$. We also say that
  the \emph{generic Torelli theorem}
  holds if $per$ has degree $1$ onto its image. The
  \emph{Schottky problem} asks for a description of the image of
  the period map.
  
  As explained on p. 228 of \cite{G},
  if the generic local Torelli theorem holds
  and if it can be proved that
  a variety $X$
  can be recovered from knowledge of the period point
  $per(X)$ and the subspace $per_*(T_\sM(X))$ of the tangent space
  $T_\sP(X)$,
  then the generic Torelli theorem also holds.

  In this paper we consider the problem for elliptic surfaces
  $f:X\to C$
  with no multiple fibres
  (such surfaces we shall call \emph{simple})
  and show that the situation is closely parallel to that for curves,
  as follows.

    Suppose that
    the geometric genus of $X$ is $h$ and its irregularity $q$,
    that $b_1(X)$ is even (so that,
    by a result of Miyaoka \cite{Mi},
    $X$ is K{\"a}hler),
    that $10 h>12(q-1)$
    and that $h\ge q+3$. Assume also that
    $X$ is \emph{general}, in a sense to be
    made precise later. Then we prove
    the following results, the first two
    of which are well known tautologies.

    For $a\in C$ we let
  $E_a$ denote the fibre $f^{-1}(a)$
  and $\omega_a^\vee$ the line in
  $H^{2,0}(X)^\vee$ whose kernel
  is $H^0(X,\Omega^2_X(-E_a))$.

  Assume that $10 h>12(q-1)$.
    
    \begin{enumerate}
  \item There is a classifying
    morphism $\phi=\phi_f:C\to\bsm_1$, where
    $\bsm_1$ is the stack of stable curves of genus $1$.
    Set $Z=\Ram_\phi\subset C$, the ramification
t    divisor. (Up to noise which is removed by
    the language of stacks, this is the locus where
    the derivative of the $j$-invariant vanishes.)
    \item
    If also $X$ is algebraic, then the tangent
    space at the point $X$ to the stack of
    algebraic elliptic surfaces is naturally isomorphic
    to an invertible sheaf on $Z$.
  \item
    Every choice of a point $a$
    in $Z$ and of a local co-ordinate on $C$ at $a$
    defines a $1$-parameter variation
    of $X$. This is based on the construction
    of the version of \emph{Schiffer variations}
    that is described on p. 443 of \cite{Ga}.
  \item
        The derivative of the period map
    of this variation,
    which is a linear map
    $\bnabla_a:H^{2,0}(X)\to H^{1,1}_{prim}(X),$
    is of rank $1$.
  \item
    There is a meromorphic $2$-form
    $\eta_a\in H^0(X,\Omega^2_X(2E_a))$
    of the second kind (that is,
    the residue of $\eta_a$ along $E_a$ vanishes)
    such that $\bnabla_a=\omega_a^\vee\otimes[\eta_a]$.
    
  \item Assume also that $h\ge q+3$. 
    Then the canonical model of $X$ is a copy of
    $C$ embedded as a curve of degree $h+q-1$
    in a projective space $\P(H^{2,0}(X)^\vee)=\P^{h-1}$, and the
    set $Z$ can be recovered, as a finite point set in
    $\P^{h-1}$, from
    the finite
    subset $\left\{\omega_a^\vee\otimes[\eta_a]\right\}_{a\in Z}$
    of $\P((H^{2,0}(X)^\vee\otimes H^{1,1}_{prim}(X))^\vee)$.
    Indeed, we exploit this set of $N$ points
    in projective space
    as an analogue
    of the theta divisor on the Jacobian of a curve.
  \item Assume that $h\ge q+3$ and that $4h>5(q-1)$.
    Then the curve $C$ can be recovered from $Z$
    in $\P^{h-1}$ via quadratic interpolation.
  \item Given $C$ and $Z$,
    we then prove a generic Torelli theorem for
    Jacobian elliptic surfaces.
  \end{enumerate}

  \begin{remark} It is clear that some of these
    constructions can be still be made when the phrase
    ``elliptic curve'' is replaced by ``Calabi--Yau variety
    whose compactified moduli stack is a smooth $1$-dimensional
    Deligne--Mumford stack whose first Chern class is positive.''
    \end{remark}

    An essential difference between the case
    of curves and that of elliptic surfaces,
  however, is that for curves these variations
  arise for any point $a$ on $C$
  while for surfaces
  they only arise for points of the ramification
  divisor $Z$. Indeed, for other points $x$ of $C$
  there is no meromorphic
  $2$-form of the second kind with double
  poles along $E_x$. (I am grateful to Richard Thomas
  for explaining this to me.)
  
  We now give some more details.

  \begin{definition}  An elliptic surface is
    \emph{Jacobian} if it has a specified
    section.
  \end{definition}
  Jacobian implies simple but not conversely.
  
 In this paper the things of primary concern are
 the stacks $\sse$ and $\tsje$ of 
 simple and Jacobian elliptic surfaces $f:X\to C$
 that are smooth and
 relatively minimal.
 We also consider the stack $\sje$ whose objects
 are the relative canonical models of surfaces
 in $\tsje$; given $f:X\to C$ in $\tsje$
 the relative canonical model
 is obtained by contracting all vertical
 $(-2)$-curves in $\X$ that are disjoint
 from the given section. There is an obvious morphism
 $\tsje\to\sje$ that is a bijection on geometric
 points. At the level of miniversal deformation
 spaces, this morphism can be described
 by taking the geometric quotient by an action
 of the relevant Weyl group, as was shown by Artin
 and Brieskorn.

We shall say that a surface
$X$ in $\sse$ is \emph{general} if its $j$-invariant
$C\to\P^1_j$ is non-constant and its singular fibres are all of type
$I_1$. We let $\sse^{gen}$ and $\sje^{gen}$
denote the stacks of general simple surfaces
and general Jacobian surfaces; these are open substacks
of $\sse$ and $\tsje$.
Note that $\sje^{gen}$ maps
isomorphically to its image in $\sje$ so is also
naturally an open substack of $\sje$.
Then to give a point in $\sse^{gen}$ is equivalent to giving a
classifying morphism $F:C\to\bsM_1$
that is non-constant and unramified over $j=\infty$.
The stack $\bsM_1$ is not the same as the stack $\bsell$
of stable generalized elliptic curves; these
stacks will be discussed in more detail in
Section \ref{2}. Giving a point in
$\sje^{gen}$ is equivalent to giving a 
morphism $\phi:C\to\bsell$ 
that is non-constant over the $j$-line
and unramified over $j=\infty$.

Assume that $f:X\to C$ is general in $\sje$.
Let $Z$ denote the ramification divisor in $C$ of $\phi$.
Say $h=p_g(X)$ and $q=h^1(\sO_X)$, so that
$q$ is also the genus of $C$.
We shall assume throughout this paper that
\begin{equation}\label{(**)} {10 h> 12 (q-1)} \
  {\textrm{and}}\ {h  >0.}
\end{equation}
These assumptions ensure that
$\deg\phi^*T_{\bsell}> 2q-2$, which in turns
ensures the vanishing
of certain obstruction spaces. From Section
\ref{recover} onwards we shall make the stronger
assumptions that
\begin{equation}\label{(***)}
  {8h > 10(q-1)} \ {\textrm{and}}\
  {h>q+2.}
\end{equation}
These assumptions make it possible
to apply theorems of Mumford and Saint-Donat
about the defining
equations of linearly normal projective curves.

Write $N=10h+8(1-q)$. Then, as is well known, $\sse$ 
is smooth at the point corresponding to $X$ and 
$$\deg Z=N,\ \dim\sse=N+h,\   h^{1,1}(X)=N+2.$$
Then we shall prove effective forms of both a generic local Torelli theorem
and a generic Torelli theorem for the weight $2$ Hodge structure on $X$, 
in the following sense.

If $X$ is a surface in $\sje$ with specified section
$\sigma$ and fibre $\xi$ then $H^2_{prim}(X)$ and
$H^{1,1}_{prim}(X)$ will denote the orthogonal complement
$\langle \s,\xi\rangle^\perp$.
If $X$ is in $\sse$
but is not necessarily algebraic then
 $H^2_{prim}(X)$ and
 $H^{1,1}_{prim}(X)$ will denote
 $\xi^\perp/\Z\xi$; these two definitions are
 equivalent for $X\in\sje$. Observe that $h^{1,1}_{prim}(X)=N$.
 In fact $\dim\sje=N$ also, so that
 $\dim\sje =h^{1,1}_{prim}(X)$. We shall use this coincidence
 in Section \ref{section linearize}
 to enhance the local structure of the derivative
 of the period map.

From the description of $\sje^{gen}$ as the stack that parametrizes
those non-constant morphisms from curves
to $\bsell$
that are unramified over $j=\infty$
we shall prove the following theorem, which appears
as Theorem \ref{deriv formula}. It is the main
result of the paper; everything else follows from it.

\begin{theorem}\label{main}
  Fix a surface $f:X\to C$ that is a point $X$ of $\sje^{gen}$
  and corresponds to $\phi:C\to\bsell$.
  For $P\in C$ put $E_P=f^{-1}(P)$ and let
  $Z$ denote  the ramification divisor $\Ram_\phi$.
  
  \part[i] Given a point $a$ of $Z$
  there is a tangent line $\C v_a$
  to $\sje$ at the point $X$.
  \part[ii] There is a meromorphic $2$-form
  $\eta_a\in H^0(X,\Omega^2_X(2E_a))$
  of the second kind (that is, its residue
  $\Res_{E_a}\eta_a$ vanishes).
  \part[iii] The corresponding map
  $$\bnabla_{v_a}:H^{2.0}(X)\to H^{1,1}_{prim}(X)$$
  is of rank $1$. Its kernel is the space
  $H^0(X,\Omega^2_X(-E_a))$ of $2$-forms
  that vanish along $E_a$ and its image
  is the line generated by the class
  $[\eta_a]$ of $\eta_a$, modulo $\xi$.
  \noproof
\end{theorem}

In order to prove this theorem
we shall use Schiffer variations
to construct, for each point $a\in \Ram_\phi$
where the ramification index
of $\phi:C\to\bsell$ is $m$
(so that $a$ is of multiplicity $m-1$
in $Z$) an $(m-1)$-parameter
deformation $\sC\to\Delta^{m-1}$
of $C$ whose derivative
can be calculated.
So, when $X$ is a point of $\sje^{gen}$,
we have a detailed description of an $N$-dimensional
subspace of the tangent space to moduli inside
the tangent space to the period domain
as the subspace spanned by certain explicit
tensors of rank $1$.
(Masa-Hiko Saito \cite{S} has proved the local Torelli
theorem for simple elliptic surfaces with non-constant $j$-invariant
and for many surfaces with constant $j$-invariant.
We shall extend his result slightly; see Theorem \ref{2.11}
below.)

It is a matter of linear algebra to recover
$Z$ as a subset of the projective space
$\P^{h-1}$ in which $C$ is embedded as the canonical model of $X$,
under the assumption that $Z$ is reduced.
We then use a theorem of Mumford \cite{Mu}
and Saint-Donat \cite{SD}, to the effect
that linearly normal curves of genus $q$
and degree at least $2q+2$
are intersections of quadrics,
to show that $C$ is determined by
quadratic interpolation through $Z$.
We go on to prove
that from the pair $(C,Z)$ we can recover
the classifying morphism
$\phi:C\to\bsell$,
modulo
the automorphism group
$\GG_m$ of $\bsell$ provided
that $\phi$ is generic.
This recovery of $C$ and $\phi$ from the period data
we regard as an effective theorem.
It shows that any failure of generic
Torelli for Jacobian surfaces
can be detected in a pencil
that is the closure of the $\GG_m$-orbit
thorough a generic point of $\sje^{gen}$.
(The fact that the automorphisms of $\bsell$
obstruct a direct deduction of generic Torelli
from knowledge of $C$ and $Z$
was observed by Cox and Donagi \cite{CD}.)

Once we know that the base curve $C$ is determined
by Hodge-theoretical data of weight $2$
we go on to
prove the generic Torelli theorem for
Jacobian elliptic surfaces
via ideas similar to those used by Chakiris \cite{C1} \cite{C2}
to prove generic Torelli when $C=\P^1$,
but reinforced by the
Minimal Model Program. 

There is also some further structure on
the period map for $\sje^{gen}$: the relevant
vector bundles and homomorphisms
between them can be described in terms of
line bundles on the universal
ramification divisor
$\sZ^{gen}$ over $\sje^{gen}$
of the universal classifying
morphism to $\bsell$. This
can be seen as a local solution to the Schottky problem.
The
details are stated in Theorem \ref{linearize}.

We also give a variational form of a
partial solution to the global Schottky problem. 
  
  If $\sX$ is a Deligne--Mumford stack then
  $[\sX]$ will denote its geometric quotient.
\end{section}

\begin{acknowledgments}
  I am very grateful to Phillip Griffiths, Ian Grojnowski, Richard Thomas and Tony Scholl
  for some valuable discussion.
\end{acknowledgments}
\begin{section}{Preliminaries on stacks and tangent spaces}\label{2}
  Everything in the next two sections is well known;
  if it is not due to either Kas \cite{Ka} or Kodaira \cite{Ko}
  then it is folklore.
  
  The stack $\bsEll$ is the Deligne--Mumford stack over $\C$
  of \emph{stable generalized elliptic curves};
  that is, an $S$-point of $\bsEll$ is a flat projective morphism $Y\to S$
  with a section $S_0$ contained in the relatively smooth locus of $Y\to S$
  and whose geometric fibres are reduced and irreducible
  nodal curves of arithmetic genus $1$.
  Such a curve is then, locally on $S$, a plane cubic with affine equation
  $$y^2=4x^3-g_4x-g_6,$$
  where $g_4$ and $g_6$ are not both zero,
  so that $\bsEll$ is the quotient stack $\P(4,6)=(\A^2-\{0\})/\GG_m$, where
  $\GG_m$ acts on $\A^2$ with weights $4,6$. Note that
  $\GG_m$ acts on $\A^2$ via a homomorphism $\GG_m\to \GG_m^2$
  (whose kernel is $\mu_2$)
  and the standard
  action of $\GG_m^2$ on $\A^2$,
  so that there is a residual action of $\GG_m$ on $\bsEll$.
  This exhibits $\GG_m$ as the full automorphism group of $\bsEll$.
  Cf. Theorem 8.1 of \cite {BN} and the calculation there on p. 139.

  The geometric quotient $[\bsEll]$ of $\bsEll$ is the compactified $j$-line $\P^1_j$;
  if $\rho:\bsEll\to\P^1_j$ is the quotient morphism then
  the automorphism group of each fibre of $\rho$ is $\Z/2$, except over $j=j_6=0$,
  where it is $\Z/6$, and over $j=j_4=1728$, where it is $\Z/4$.
  So $\deg\rho=1/2$. As is well known, it is possible to write
  down a generalised elliptic curve  over the open locus $U$ of $\P^1_j$
  defined by $j\ne 0,1728$, so that
  there is a section of $\rho$ over $U$. Moreover,
  $\rho^{-1}(U)$ is isomorphic to $U\times B(\Z/2)$,
  but there is no global section of $\rho$.

There are two obvious line bundles on $\bsEll$:
the bundle $M$ of modular forms of weight $1$, which is identified with
the conormal bundle of the zero-section of the universal stable
generalized elliptic curve, and the tangent bundle
$T_{\bsEll}$.

\begin{lemma} 
  $T_{\bsEll}\cong M^{\otimes 10}$, $\deg M=1/24$ and
  $\deg T_{\bsEll}=5/12$.

  \begin{proof} Quite generally the Picard group of
    $\P(a,b)$ is generated by $\sO(1)$, which has degree
    $1/ab$, and $T_{\P(a,b)}$ is isomorphic to
    $\sO(a+b)$, which then has degree $(a+b)/ab$.
  \end{proof}
\end{lemma}

The objects of
the stack $\bsM_1$ are
\emph{stable curves of genus $1$};
the geometric fibres are isomorphic to stable generalised
elliptic curves, but no section is given. This is an Artin stack,
but not Deligne--Mumford. Indeed, the word ``stable''
in this context
is an abuse of language, but I am optimistic that it will
cause no confusion.

Let $\sC\to\bsM_1$ and $\sE\to\bsEll$ denote
the universal objects and let $\sG\to\bsEll$ denote
the N{\'e}ron model of $\sE\to\bsEll$, so that
$\sG$ is the open substack of $\sE$ obtained
by deleting the singular point of the fibre over $j=\infty$.

The next result is due to Altman and Kleiman \cite{AK},
although reformulated here in the language of stacks.
We have chosen to include a slightly different proof that emphasizes
automorphism groups rather than Picard varieties
so that the relevant classifying stacks enter more easily.

\begin{theorem}
  There is a morphism
  $\pi:\bsM_1\to\bsEll$ via which $\bsM_1$ is isomorphic
  to the classifying stack $B\sG$ over $\bsEll$.
  \begin{proof}
    Let $G\to \bsM_1$ denote the connected component
of the relative automorphism group scheme
of $\sC\to \bsM_1$.
So $G\to\bsM_1$ is elliptic over the open
substack $M_1$ of $\bsM_1$ defined by $j\ne\infty$
and over $j=\infty$ the fibre of $G$ is the multiplicative group $\GG_m$.

\begin{lemma} There is a unique open
embedding
$G\inj\tG$, over $\bsM_1$,
where $\tG\to \bsM_1$
is a stable generalised elliptic curve and $G$ is its
relative smooth locus.

\begin{proof} To construct $\tG$
  we need to patch the puncture
  of $G\to \bsM_1$ that lies over $j=\infty$.
  In a suitable \nbd $S$ of the locus $j=\infty$ in $\bsM_1$
  the process of
  patching the puncture
  is a matter of ``reversing the process of deleting a
  closed point
  from a normal $2$-dimensional analytic space (or scheme)'',
  so the patch is unique
  if it exists. Therefore it is enough to exhibit the patch
  locally on $\bsM_1$, in the \nbd $S$.
  
  Put $\sC_S=\sC\times_{\bsM_1}S$.
Since $S$ is local and $\sC_S\to S$ is
generically smooth, there is a section over $S$ of
$\sC_S\to S$ that is contained in the relative
smooth locus.
Use this section to put the structure of a stable generalised elliptic
curve on $\sC_S\to S$. This
structure on $\sC_S\to S$ provides the patch for
$G_S=G\times_{\bsM_1}S\to S$ and yields $\tG\to S$;
the lemma is proved.
\end{proof}
\end{lemma}

Sending $\sC\to \bsM_1$ to $\tG\to \bsM_1$
defines a morphism $\pi: \bsM_1\to \bsEll$.

A \emph{gerbe} is a morphism of stacks that is
  locally surjective on both objects and morphisms.
  A \emph{neutral gerbe} is a gerbe with a section.

\begin{lemma}
$\pi$ is a neutral gerbe.
\begin{proof} 
  It  is enough to
  show that
  
\noindent (1) $\pi$ has a section (then it is certainly locally surjective on objects) and

\noindent (2) $\pi$ is locally surjective on morphisms.

For (1), 
use the forgetful morphism
$\bsEll\to\bsM_1$ to get a section of $\pi$.

(2) is equally clear: a morphism of stable generalised elliptic curves 
is, in particular, a morphism of the underlying stable curves.
\end{proof}
\end{lemma}

Finally, suppose that $\sX\to \sY$ is a neutral gerbe with section
$s:\sY\to \sX$. The stabiliser group scheme
is a group scheme
$\sH'\to \sX$; define $\sH\to \sY$ to be the pull back of $\sH'$
via $s$. If $\sH\to \sY$
is flat then \cite{L--MB} there is an isomorphism
$\sX\to B\sH$.
Since the N{\'e}ron model $\sG\to\bsEll$ is flat, the theorem is proved.
\end{proof}
\end{theorem}

Note that $G\to\bsM_1$ is isomorphic to the pull back
under $\pi$ of the N{\'e}ron model $\sG\to\bsEll$.

As already remarked,
a general simple surface $f:X\to C$ determines, and is determined by,
a morphism $F=F_f:C\to\bsM_1$ that is unramified over
$j=\infty$. Let
$\phi=\phi_f=\pi\circ F:C\to\bsEll$ denote the composite,
so that the induced Jacobian elliptic surface is the compactified
relative automorphism group scheme of $f:X\to C$.

For example, if $f:X\to \P^1$ is a primary Hopf surface
then $\phi_f$ is constant: the relative automorphism group
scheme is a constant relative
group scheme $E\times\P^1\to\P^1$.

\begin{lemma} If $X\to C$ has non-constant
  $j$-invariant then the irregularity $q$
of $X$ equals the geometric genus of $C$.
\begin{proof}
  This is well known, and easy.
  \end{proof}
\end{lemma}

\begin{lemma} Suppose that $f:X\to C$ is a general
  simple surface.
  
  \part[i] $\deg\phi= 2 c_2(X)=24\chi(X,\sO_X)$.
  \part[ii] $\deg\phi^*M=\chi(X,\sO_X)$.
  \part[iii] $\deg\phi^*T_{\bsEll}=10\chi(X,\sO_X)$.
\begin{proof} Since $c_2(X)$ equals the total number of nodes
  in the singular fibres of $X\to C$,
  \DHrefpart{i} follows from consideration
of the inverse image of the locus $j=\infty$, the fact that $\deg\rho=1/2$
and Noether's formula. The rest follows immediately.
\end{proof}
\end{lemma}

We next consider various tangent spaces.

Regard $B\sG$ as the quotient of $\bsEll$ by $\sG$.
Since $\bsm_1$ is isomorphic to $B\sG$, locally
on $\bsell$,
the tangent complex
$T^\bullet_{\bsm_1}$
is a $2$-term complex, which is obtained by descending
a $2$-term complex on $\bsEll$. In degree $0$ this complex
is $T_{\bsEll}$, in degree $-1$ it is the adjoint bundle $\Ad \sG$
and the differential is the derivative of the action
of $\sG$ on $\bsEll$.

Since this action is trivial,
the differential in $T^\bullet_{\bsm_1}$
is zero. Moreover, since $\sG$ has no characters
(it is generically an elliptic curve) it follows that
$T^\bullet_{\bsm_1}$ is quasi-isomorphic to the pull back of
the complex $\Ad \sG[1]\oplus T_{\bsEll}[0]$ on $\bsEll$.

Note that $(\Ad\sG)^\vee$ is exactly the line bundle $M$
of modular forms of weight $1$.

Now fix a point $X$ of $\sse$.
That is, we fix a general simple elliptic surface $f:X\to C$.
This equals the datum of a morphism
$F:C\to \bsM_1$. Set $\phi=\pi\circ F:C\to\bsEll$. Then
$F^*T^\bullet_{\bsm_1}$ is quasi-isomorphic to
$\phi^*\Ad \sG[1]\oplus \phi^*T_{\bsEll}$. 
There is a distinguished triangle
$$T_C\to F^*T^\bullet_{\bsm_1}\to K^\bullet$$
where $K^\bullet$ is a $2$-term complex
of coherent sheaves on $C$,
$K^{-1}=\phi^*\Ad \sG$, $K^0$ is the skyscraper sheaf
$\coker (T_C\to\phi^*T_{\bsEll})$
and the differential in $K^\bullet$ is zero.

\begin{proposition}\label{2.8} The tangent space $T_{\sse}(X)$ is naturally
  isomorphic to the hypercohomology group
  $\H^0(C,K^\bullet)$ and the obstructions
  to the smoothness of $\sse$ at $X$ lie in
  $\H^1(C,K^\bullet)$.
  \begin{proof} This follows from the identification of the points of
    $\sse$ with morphisms from curves to $\bsm_1$.
    In this latter context the result is well known.
  \end{proof}
\end{proposition}

Let $Z$ denote the ramification divisor
$Z=\Ram_\phi=\Ram_F$ on $C$.
Then
$K^0$ is an invertible sheaf on $Z$,
so that there is a non-canonical
isomorphism $K^0\buildrel\cong\over\to \sO_Z$.
    
  The distinguished triangle just mentioned gives
  an exact sequence
\begin{eqnarray*}
{0}&\to& H^0(C,T_C) \to\H^0(C,F^*T^\bullet_{\bsm_1})\to\H^0(C,K^\bullet)\\
        &\to &{H^1(C,T_C)}  \to \H^1(C,F^*T^\bullet_{bsm_1})\to\H^1(C,K^\bullet)\to 0.
\end{eqnarray*}
\noindent Then $H^2(C,\phi^*\Ad G)=0$ and
$H^1(C,\phi^*T_{\bsEll})=0$, since
$$\deg\phi^*T_{\bsEll}=10\chi(X,\sO_X)>2q-2,$$
from the assumption (\ref{(**)}),
so that
$\H^1(C,K^\bullet)=\H^1(C,F^*T^\bullet_{\bsm_1})=0$ and there is 
an exact commutative diagram
$$\xymatrix{
  &{0}\ar[d]&{0}\ar[d]&&\\
  &{H^0(C,T_C)}\ar[r]^{=}\ar[d]&{H^0(C,T_C)}\ar[d]&&\\
  {0}\ar[r]&{H^0(C,\phi^*T_{\bsell})}\ar[r]\ar[d]&
  {\H^0(C,F^*T^\bullet_{\bsm_1})}\ar[r]\ar[d]&
    {H^1(C,\phi^*\Ad\sG)}\ar[r]\ar[d]^{=}&{0}\\
    {0}\ar[r]&{H^0(C,K^0)}\ar[r]\ar[d]&
    {\H^0(C,K^\bullet)}\ar[r]\ar[d]&
    {H^1(C,\phi^*\Ad\sG)}\ar[r]&{0}\\
    &{H^1(C,T_C)}\ar[r]^{=}\ar[d]&{H^1(C,T_C)}\ar[d]&&\\
    &{0}&{0}&&
  }$$
  in which the two middle rows are canonically split.

\begin{proposition}\label{2.4}\label{2.9}

\part[i] The stack $\sse$ is smooth at the point $X$.

\part[ii] Its dimension there is $11h+8(1-q)$.

\part[iii] The degree of the ramification divisor $\Ram_\phi=\Ram_F$ is $N$.

\part[iv] $\phi^*M\cong
f_*\omega_{X/C}\cong(R^1f_*\sO_X)^\vee$
and $\deg\phi^*M=\chi(X,\sO_X)=h+1-q$.
\begin{proof} Except for \DHrefpart{iv},
  which is well known, this follows from the preceding discussion.
\end{proof}
\end{proposition}

Suppose that $f:X\to C$ is a point of $\sje^{gen}$ and defines
$\phi:C\to\bsell$.
Set $Z=\Ram_\phi$.

\begin{proposition}\label{2.5}

  \part[i]
  There is a short exact sequence
  $$0\to H^0(C,\phi^*T_{\bsell})\to T_{\sje}(X)
  \to H^1(C,T_C)\to 0.$$
  
  \part[ii] $T_{\sje}(X)$ is
  naturally isomorphic to
  $H^0(Z,K^0)$.
  
  \part[iii] $\sje$ is smooth at the point $X$ and its dimension
  there is $N$.
  \begin{proof} This is proved in the same way
    as Propositions \ref{2.8} and \ref{2.9}.
  \end{proof}
\end{proposition}

Since $Z$ is $0$-dimensional and $K^0$
is an invertible sheaf on $Z$,
$H^0(Z,K^0)$ is non-canonically isomorphic to $H^0(Z,\sO_Z)$.

We shall usually write $\phi^*M=L=f_*\omega_{X/C}$. Then
$\phi^*T_{\bsEll}\cong L^{\otimes 10}$, so that
$$\sO_C(\Ram_\phi)\cong \phi^*T_{\bsEll}\otimes T_C^\vee
\cong L^{\otimes 10}\otimes\sO_C(K_C).$$
At this point we give a slight refinement of M.-H. Saito's
local Torelli theorem. The argument is essentially his.

\begin{theorem}\label{saito}\label{2.11}
  Suppose that $f:X\to C$ is a simple elliptic
  surface with $r$ singular fibres. Put $L=f_*\omega_{X/C}$.
  Assume
  that $r\ge\deg L+3$ and that
  $\deg L\ge 2$. Then the local Torelli theorem
      holds for $X$.
      \begin{proof} From the main result of \cite{S}
        it is enough to consider the situation where the $j$-invariant
        is constant.
        Following \cite{S} it is enough to show that the natural
      homomorphisms
      $$\mu_1: H^0(C,f_*\Omega^2_X)\otimes H^1(C,f_*\Omega^1_X)
      \to H^1(C,f_*(\Omega^1_X\otimes\Omega^2_X))\ \textrm{and}$$
      $$\mu_2: H^0(C,f_*\Omega^2_X)\otimes H^0(C,R^1f_*\Omega^1_X)
      \to H^0(C,R^1f_*(\Omega^1_X\otimes\Omega^2_X))$$
      are surjective. Recall that $K_X\sim f^*(K_C+L)$,
      so that
      $f_*(\Omega^1_X\otimes\Omega^2_X)\cong
      \sO_C(K_C+L)\otimes f_*\Omega^1_X$,
      and that $\vert K_C+L\vert$ has no base points.

      \begin{lemma}\label{bpf} If $\sF,\sG$ are coherent sheaves
        on a $1$-dimensional projective scheme $C$ over a field $k$
        and if $\sF$
        is generated by $H^0(C,\sF)$, then the natural multiplication
        $$H^0(C,\sF)\otimes_k H^1(C,\sG)\to H^1(C,\sF\otimes_{\sO_C} \sG)$$
        is surjective.
        \begin{proof} There is an exact sequence
          $$H^0(C,\sF)\otimes_k\sO_C \to \sF\to 0;$$
          tensoring this with $\sG$ gives an exact sequence
          $$H^0(C,\sF)\otimes_k\sG\to \sF\otimes \sG\to 0.$$
          Taking $H^1$ of this sequence
          gives the result.
        \end{proof}
      \end{lemma}
      In particular, taking
      $\sF=\Omega^1_C\otimes L=f_*\Omega^2_X$
      and $\sG=f_*\Omega^1_X$ shows that
      $\mu_1$ is surjective.

      Now consider $\mu_2$. Set $A=\sum_1^r a_i$, the critical
      subset of $C$. The exact sequences (4.17) and (4.18)
      of \cite{S} are
      $$0\to \Omega^1_C\to f_*\Omega^1_X
      \to\sO_C(L-A)\to 0\ \textrm{and}$$
      $$0\to\sO_C(K_C+A-L)\to R^1f_*\Omega^1_X
      \to \sO_C\oplus \sT^1\to 0,$$
      where $\sT^1$ is a skyscraper sheaf.
      The second sequence then gives
      $$0\to\sO_C(K_C+A-L+\delta)\to R^1f_*\Omega^1_X
      \to \sO_C\to 0,$$
      where $\delta\ge 0$,
      via the process of saturating a subsheaf.
      Since $\deg A>\deg L$, by assumption, this last
      sequence splits and
      $$R^1f_*\Omega^1_X=\sO_C\oplus \sO_C(K_C+A-L+\delta).$$
      From its definition, $\mu_2$ is then the direct sum
      $$\mu_2=\mu_2'\oplus\mu_2''$$
      of multiplication maps 
      $$\mu_2':H^0(C,\sO_C(K_C+L))\otimes\C
      \to H^0(C,\sO_C(K_C+L))\ \textrm{and}$$
      $$\mu_2'':H^0(C,\sO_C(K_C+L))\otimes
      H^0(C,\sO_C(K_C-L+A+\delta))
      \to H^0(C,\sO_C(2K_C+A+\delta)).$$
      The first of these is obviously surjective,
      while the surjectivity of the second follows from
      \cite{Mu}, Theorem 6, p. 52 and the facts that
      $\deg(K_C+L)\ge 2q$ and
      $\deg(K_C-L+A+\delta)\ge 2q+1$.
    \end{proof}
  \end{theorem}
  \end{section}
\begin{section}{The comparison between $\sse^{gen}$ and $\sJE^{gen}$}\label{comparison}
The morphism $\pi:\bsM_1\to\bsEll$
defines a morphism $\Pi:\sse^{gen}\to\sje^{gen}$. If we fix
a point $f:X\to C$ of $\sje^{gen}$, then $\Pi^{-1}(X)$
is identified with $H^1(C,\sH)$, where $\sH\to C$ is the
N{\'e}ron model of $X\to C$.

The sheaves $L^\vee$ and
$\sH$ are group schemes over $C$. Define $\sF=R^1f_*\Z$; this
is a constructible sheaf of $\Z$-modules on $C$ that is of generic
rank $2$.

The following results are due, in essence,  to Kodaira.
In particular, 
Proposition \ref{kod} is a variant of \cite{Ko}, p. 1341, Theorem 11.7;
there he proves only that $H^2(C,\sF)$ is finite,
but he does not assume that $f:X\to C$ is general.

\begin{lemma} There is a short exact sequence
  $$0\to\sF\to L^\vee\to\sH\to 0$$
  of sheaves of commutative groups on $C$.
  \begin{proof}
    This follows from the exponential exact sequence on $X$
    and the identification
    $\sH=\ker(R^1f_*\sO_X^*\to R^2f_*\Z)$.
  \end{proof}
\end{lemma}

\begin{proposition}\label{kod}
  $H^2(C,\sF)=0$
  and the homomorphism
  $H^1(C,L^\vee)\to H^1(C,\sH)$ is surjective.
  \begin{proof}
    It is enough to show that $H^2(C,\sF)=0$.

    For any ring $A$, there is a Leray spectral sequence
    $$E^{pq}_{2,A}=H^p(C,R^qf_*A)\Rightarrow H^{p+q}(X,A).$$
    Since the fibres of $f$ are irreducible curves
    and $X\to C$ has a section,
    the hypotheses of
    Th{\'e}or{\`e}me 1.1 of
    \cite{SGA7II} XVIII are satisfied, so that this degenerates
    at $E_2$.
    Take $A=\Z$; then
    there is a short exact
    sequence
    $$0\to H^2(C,\sF)\to H^3(X,\Z)\to H^1(C,\Z)\to 0.$$
    So $H^2(C,\sF)$ is identified with the
    torsion subgroup $\Tors H^3(X,\Z)$ of $H^3(X,\Z)$.

    Suppose that $\ell$ is a prime dividing the order of
    $\Tors H^3(X,\Z)$. Since $H^4(X,A)$ is isomorphic to $A$, it follows
    from taking cohomology of the short exact sequence
    $$0\to\Z\buildrel{.\ell}\over\to\Z\to\Z/\ell\to 0$$
    that $H^3(X,\Z/\ell)\cong H^3(X,\Z)\otimes\Z/\ell.$
    Also, Poincar{\'e} duality gives an isomorphism
    $ H^3(X,\Z/\ell)\to H^1(X,\Z/\ell)^\vee,$ where ${}^\vee$
    denotes the dual $\Z/\ell$-vector space.

    The spectral sequence $E^{pq}_{r,\Z/\ell}$ shows that
    $\beta:H^1(C,\Z/\ell)\to H^1(X,\Z/\ell)$ is injective, so if it
    is not surjective then the map
    $i^*: H^1(X,\Z/\ell)\to H^1(\s,\Z/\ell)$
    induced by the inclusion $i:\s\inj X$
    of the zero section is not
    injective. Then there is an {\'e}tale $\Z/\ell$-cover
    $\a:\tX\to X$ that is split over $\s$. So $\tX\to C$ is
    elliptic and has a section $\ts$ such that
    $N_{\ts,\tX}\cong N_{\s/X}$. However,
    $$\deg N_{\ts,\tX}
    =-\chi(\tX,\sO_{\tX})=
    -\ell\chi(X,\sO_X)=\ell\deg N_{\s/X}.$$
    So $\beta$ is an isomorphism, so that
    $H^3(X,\Z/\ell)\cong(\Z/\ell)^{2q}$
    and therefore $H^3(X,\Z)$ is torsion-free.
  \end{proof}
\end{proposition}

Let $\sse_{h,q}$ denote the substack of $\sse$
that consists of surfaces whose geometric genus is $h$
and whose irregularity is $q$. This is a union
of connected components of $\sse$.

The next result is an immediate corollary
of Proposition \ref{kod}.

\begin{corollary} \part[i] The closed substack $\Pi^{-1}(X)$ is
  irreducible.
  \part[ii] $\sse_{h,q}$ is
  irreducible.
  \noproof
\end{corollary}

Now suppose that $Y\in\Pi^{-1}(X)$. According to
\cite{Ko}, p. 1338, Theorem 11.5, $Y$ is algebraic
if and only if it defines a torsion element of
$H^1(C,\sH)$.

\begin{proposition} The algebraic surfaces are dense in $\Pi^{-1}(X)$.
  \begin{proof} We must show that the image of
    $H^1(C,R^1f_*\Z)\otimes\Q$ in $H^1(C,R^1f_*\sO_X)$
    is dense.

    Let $\xi^\perp$ denote the orthogonal complement
    of $\xi$ in $H^2(X,\Z)$.
    Then, via
    the Leray spectral sequence, there is a commutative square
    $$\xymatrix{
    {H^1(C,R^1f_*\Z)\otimes\Q}\ar[r]^{\a}\ar[d]_{\cong}&{H^1(C,R^1f_*\sO_X)}\ar[d]^{\cong}\\
    {(\xi^\perp/\Z\xi)\otimes\Q}\ar[r]^{\beta} &{H^2(X,\sO_X)}
    }
    $$
    where the vertical arrows are isomorphisms.
    Then $\beta$ has dense image, from the K{\"a}hler property of $X$,
    and the proposition is proved.
  \end{proof}
\end{proposition}
\end{section}
\begin{section}{Schiffer variations for elliptic surfaces
    and the derivative of the period map}\label{schiffer}
  Fix a curve $C$, a point $a\in C$, a local co-ordinate $z$ on $C$ at
  $a$, an integer $e\ge 2$
  and a real number $\delta$ such that $0<\delta \ll 1$
  and the region
  $\vert z\vert <\delta$ is an open disc $D=D_z$ in $C$.
  We begin by recalling, from \cite{Ga} p. 443, the construction of
  certain variations that we shall call
  \emph{Schiffer variations}: these
  are families $\pi:\sC\to\Delta^{e-1}=\Delta^{e-1}_{\ut}$
  over $(e-1)$-dimensional polydiscs with co-ordinates
  $\ut=(t_2,\ldots ,t_{e})$ whose closed fibre
  $\pi^{-1}(\u0)$ is $C$. They are constructed as follows.

  In $\C^{e-1}_{\ut}$ take the polydisc $\Delta^{e-1}_{\ut}$ defined
  by
  $\vert t_j\vert <\delta^{4e}$ for all $j$. In $C\times
  \Delta^{e-1}_{\ut}$
  take the complement $U$ of the closed subset
  $\bigcup_j(\vert z\vert^{4e}\le \vert t_j\vert).$
  So $U$ is a thickening of the punctured curve
  $C-\{a\}$.
  Then put $G=U\cap (\vert z\vert <\delta)$.

  In $\C^2_{z,w}\times\Delta^{e-1}_{\ut}$
  take the subset $F$ defined by
  $$
  \vert t_j\vert < \vert z\vert^{4e} \ \forall j,\
  \vert z\vert <\delta,\
  \vert w-z\vert < \vert z\vert ^{2e}\
  \textrm{and}\
  w^e+e\sum_0^{e-2} t_{e-j}w^{j}=z^e.
  $$
  Then there are projections $p:F\to U$ and
  $q:F\to \C_w\times\Delta^{e-1}_{\ut}$.

\begin{lemma}
  $p$ and $q$ are
  open embeddings.
  \begin{proof} It is enough to show that they are
    unramified and separate points.

    The ramification locus of $p$ is defined by
    $$
    w^{e-1}+\sum jt_{e-j}w^{j-1}=0.
    $$
    Since $\vert t_{e-j}\vert <\vert z\vert^{4e}$ this gives
    $\vert w\vert <\vert z\vert^4.$ Then
    $\vert z\vert-\vert z\vert^{2e}<\vert z\vert^4$,
    so that $z=0$ and $p$ is unramified.

    To check the separation of points, suppose that
    $$w_\a^e+e\sum t_{e-j}w_\a^j-z^e=0$$
    for $\a =1,2$ and $w_1\ne w_2$. Then
    $$
    \prod_{r=1}^{e-1}(w_1-\z_e^r w_2)=e\sum_{j=1}^{e-2}
    t_{e-j}\prod_{s=1}^{j-1}(w_1-\z_j^sw_2),
    $$
    where, for any integer $n$, $\z_n$ is a primitive $n$'th root of unity.
    Since $\vert w_1-z\vert<\vert z\vert^{2e}$ and
    $\vert \z_j^s w_2-\z_j^s z\vert<\vert z\vert ^{2e}$
    it follows that
    $$
    \vert w_1-\z_j^s w_2\vert <
    2\vert z\vert^{2e}+\vert z-\z_j^sz\vert
    =2\vert z\vert^{2e}+\vert 1-\z_j^s\vert\vert z\vert<
    4\vert z\vert.
    $$
    So
    $$
    \left\vert e\sum_{j=1}^{e-2} t_{e-j}\prod_{s=1}^{j-1}(w_1-\z_j^sw_2)\right\vert
    \le 4e(e-2)\vert z\vert^{4e+1}.
    $$
    On the other hand
    $\vert w_1-z\vert<\vert z\vert^{2e}$ and
    $\vert\z_e^r-\z_e^rz\vert<\vert z\vert^{2e},$
    so that
    $$
    \vert w_1-\z_e^r w_2\vert\ge \vert z-\z_e^rz\vert-
    2\vert z\vert^{2e}
    \ge \lambda\vert z\vert,
    $$
    where $\lambda=\vert\sin(2\pi/e)\vert/2.$
    So
    $$\left\vert \prod_{r=1}^{e-1}(w_1-\z_e^rw_2)\right\vert
    \ge\lambda^{e-1}\vert z\vert^{e-1}$$
    and therefore
    $$\lambda^{e-1}\le 4e(e-2)\vert z\vert^{3e+2}.$$
    This
    contradiction proves the result for $p$.

    The argument for $q$ is similar but easier,
    so we omit it.
  \end{proof}
\end{lemma}

\begin{lemma} The image of $p$ is $G$.
  \begin{proof} This is a matter of showing that,
    given $(z,\ut)\in G$, we can solve the equation
    $w^e+e\sum t_{e-j}w^j-z^e=0$
    with $\vert w-z\vert<\vert z\vert^{2e}$.
    Without the inequality there are $e$ solutions;
    if $\vert w-z\vert\ge \vert z\vert^{2e}$ for all
    of them then we get a contradiction to
    $z^e-w^e=e\sum t_{e-j}w^j$ and the inequalities
    satisfied by the $t_{e-j}$.
  \end{proof}
\end{lemma}

For each $\ut\in\Delta^{e-1}_{\ut}$ the
intersection $q(F)\cap (\C_w\times\{\ut\})$ 
is an annulus
$A_{\ut}$ in $\C_w\times\{\ut\}$ that surrounds zero.
Let $R_{\ut}$ denote the open simply connected
region in $\C_w\times\{\ut\}$ that contains $0$ and has
the same outer boundary as $A_{\ut}$, and put
$H=\sqcup_{\ut}R_{\ut}$.
So $H$ is a tubular \nbd of $\{0\}\times\Delta^{e-1}_{\ut}$
in $\C_w\times\Delta^{e-1}_{\ut}$ and $q(F)\subseteq H$.

Define $\sC$ to be the result of glueing $U$ and $H$ together
along $F$ via the maps $p$ and $q$. This is Hausdorff and, after
shrinking $\Delta^{e-1}_{\ut}$ if necessary, the morphism
$\pi:\sC\to\Delta^{e-1}_{\ut}$ is proper and is the morphism
that we sought. 
This is sometimes expressed by saying that
$\sC_{\ut}$ is constructed from
$C$ by deleting a small $z$-disc around
$a$ and glueing in a $w$-disc, where
$w$ is defined implicitly by
$w^e+e\sum_0^{e-2} t_{e-j}w^{j}=z^e.$

Until after Theorem \ref{linearize}
we fix a point
$f:X\to C$ of $\sje^{gen}$
and a point $a$ in the ramification divisor $Z=\Ram_\phi$
of the classifying morphism $\phi:C\to \bsEll$.
In particular, $\phi$ is unramified over $j=\infty$.
Put $E_a=f^{-1}(a)$ and
denote by $e=e(a)$ the ramification index at $a$ of $\phi$,
so that $a$ has multiplicity $e-1$ in $Z$. Assume
that the disc $D$ is sufficiently small,
so that it contains no other points of $Z$.
We use the
Schiffer variations of $C$
that we have just described
to construct variations of $X$.

\begin{proposition}\label{mod t2}
  For some choice of local co-ordinate $z$ on $C$ at $a$
  the morphism
  $\phi:C\to\bsell$ lifts to a morphism
  $\Phi:\sC\to\bsell$ in such a way that the restriction
  of $\Phi$ to $U$
  factors through the projection $U\to C-\{a\}$.
  \begin{proof} Given a local co-ordinate $s$
    on $\bsell$ we have $\phi^*s=z^e$
    for some local co-ordinate $z$ on $C$. Then we
    define $\Phi$ on $H$ by
    $$\Phi^*s=w^e+e\sum_{j=0}^{e-2}t_{e-j}w^{j}$$
    and on $U$ we define $\Phi$
    by composing $\phi$ with the projection $U\to C-\{a\}$.
  \end{proof}
\end{proposition}

Any Jacobian deformation of the surface $X$ determines,
for each point $a\in Z$, a deformation
of the finite scheme $V_{e(a)}=\Sp\C[t]/(t^{e(a)})$,
so that there is a morphism of local analytic deformation
spaces $\Psi:\Def_X\to \prod_{a\in Z}\Def_{V_{e(a)}}.$
Recall that $\Def_{V_{e(a)}}$ is smooth of dimension $e(a)-1$.

\begin{proposition}
  \part[i] $\Psi$ is a local analytic isomorphism.
  \part[ii] The universal ramification divisor
  $\sZ$ is smooth over $\C$.
  \begin{proof} 
    \DHrefpart{i} is a an immediate consequence of the formula
    used to define $\Phi^*s$, and
    \DHrefpart{ii} follows at once.
  \end{proof}
\end{proposition}

So we have morphisms $\sX\buildrel{F}\over\to\sC\to\Delta^{e-1}_{\ut}$
and $F:\sX\to\sC$ is a family, parametrized by $\Delta^{e-1}_{\ut}$,
of Jacobian elliptic surface whose fibre over
$0\in\Delta_{\ut}^{e-1}$ is $f:X\to C$.
   
Fix a suitable basis of $H_2(X,\Z)$,
which we identify with $H_2(\sX_{\ut},\Z)$.
Then there are holomorphic $(e+1)$-forms $\Omega^{(1)},...,\Omega^{(h)}$
on $\sX$ such that the residues
$$\omega^{(j)}(\ut)=
\Res_{\sX_{\ut}}
\left( {\begin{array}{c}{\Omega^{(j)}}\\
          {\ut}\\
          \end{array}}\right)$$
      form a normalized basis of $H^0(\sX_{\ut},\Omega^2_{\sX_{\ut}})$
      for all $\ut$.
In particular, there are $2$-cycles $A_1,...,A_h$ on $\sX_{\ut}$
such that $\int_{A_i}\omega^{(j)}(\ut)=\delta_i^j$, the Kronecker delta.

Since the line bundle
$\Omega^{e+1}_{\sX}$
pulls back from a line bundle on $\sC$,
we can expand $\Omega^{(j)}$ as
$$\Omega^{(j)}=(-1)^{e-1}\sum_{p\ge0,\uq\ge\u0} b_{p,\uq}^{(j)}
w^p\ut^{\uq}\, dw\wedge d{\ut}\wedge dv,$$
where $v$ is a fibre co-ordinate.

\begin{lemma}
  \part[i]\label{cocycle}
  $w=z(1-\sum_{i=0}^{e-2} t_{e-i}z^{i-e})
  =z(1-\sum_{i=2}^et_iz^{-i})$
  modulo $(\ut)^2$.
  \part[ii] $dw\wedge d{\ut}\wedge dv=
  \left(1+\sum_{i=2}^e(i-1)t_{i}z^{-i}\right)\,
  dz\wedge d{\ut}\wedge dv$
  modulo $(\ut)^2$.
  \begin{proof} Immediate from
    $w^e+e\sum t_{e-i}w^i=z^e$.
  \end{proof}
\end{lemma}

\begin{remark} Lemma \ref{cocycle} shows that,
  in terms of $H^1(C,T_C)$,
  the first order deformation obtained from
  $\sC\to\Delta^{e-1}_{\ut}$ is the one that arises
  from integrating the space of
  {\v{C}}ech $1$-cocycles with values in $T_C$,
  with respect to the cover
  $\{D,\ C-\{a\}\}$ of $C$,
  that is
  generated by the vector fields
  $z^{-1}d/dz,\ldots, z^{-(e-1)}d/dz$ on $D-\{a\}$.
  However, there are different kinds of Schiffer
  variation that lead to the same space of cocycles;
  this feature is part of their strength.
\end{remark}

We substitute this into the expansion of $\Omega^{(j)}$.
Since
$$\left(1-\sum t_iz^{-i}\right)^p\left(1+\sum(i-1)t_iz^{-i}\right)=
1+\sum(i-p-1)t_iz^{-i}$$
modulo $(\ut)^2$,
we get
\begin{eqnarray}\label{expand}
  \omega^{(j)}(\ut)=\sum
  b_{p,\uq}^{(j)}
z^p
\left(1+\sum(i-p-1)t_iz^{-i}\right)\ut^{\uq}dz\wedge dv.
\end{eqnarray}
\noindent modulo $(\ut)^2$.

Next, write $\omega^{(j)}(t)=\omega^{(j)}+
\sum_{i=0}^{e-2} t_{e-i}\eta_{e-i}^{(j)}$
modulo $(\ut)^2$, where
$$\omega^{(j)}=\omega^{(j)}\vert_{\ut=\u0}\
\textrm{and}\ \eta_i^{(j)}=\frac{\partial}{\partial t_i}\omega^{(j)}\vert_{\ut=\u0}.$$
Moreover, every $2$-cycle $\g$
on $X$ that is disjoint from $E$ is identified, via a
$C^\infty$ collapsing map, with a $2$-cycle on $\sX_{\ut}$
and
\begin{eqnarray}\label{star}
  \int_\g\omega^{(j)}(t)=\int_\g\omega^{(j)}+\sum_i t_i\int_\g\eta_i^{(j)}\
  \textrm{for\ all\ such}\ \g.
\end{eqnarray}

Consider the class $[\omega^{(j)}(t)]\in H^2(\sX_{\ut},\C)$. Let
$\sE\subset\sX_{\ut}$ be a fibre. Then
$[\omega^{(j)}(t)]\in\sE^\perp$.
Moreover, from the exact sequence
$$\Z\xi=\Z[E_a]\to H^2(X,\Z)\to H^2(X-E_a,\Z)$$
and the formula (\ref{star}) it follows that $\eta_i^{(j)}$ defines a class
$$[\eta_i^{(j)}]\in \xi^\perp/\Z\xi=H^2_{prim}(X)$$
such that
$$\frac{\partial}{\partial t_i}
\left([\omega^{(j)}(t)]\right)\vert_{t=0}=[\eta_i^{(j)}]\pmod{\Z\xi}.$$
Griffiths transversality shows that in fact $[\eta_i^{(j)}]$ lies in
$Fil^1=Fil^1H^{2}_{prim}(X)$. (We let $Fil^i$ refer to the $i$'th
piece of the Hodge filtration of $H^{2}_{prim}(X)$.)

We have constructed, for each $a\in Z$ of ramification index
$e(a)=e$,
an $(e-1)$-parameter variation
$\sX \to \Delta^{e-1}$ of $X$ such that the tangent space
$T_{\Delta^{e-1}}(\u0)$ is identified with the
cyclic skyscraper sheaf $L_a$ of length $e-1$ on $C$
that is supported at $a$
and determined by $Z$.
For any $v\in L_a$
consider the derivative
$$\nabla_v: H^{2,0}(X)=Fil^2\to Fil^1.$$ 
Let
$$\bnabla_v:H^{2,0}\to H^{1,1}_{prim}(X)$$
denote the composite of $\nabla_v$ with the projection
$Fil^1\to H^{1,1}_{prim}(X).$ If $L_a$ is of length $1$
then we write $\nabla_a$ rather than $\nabla_v$.

Assume that the geometric genus $h$ of $X$ is not zero.
Also, given $\omega\in H^0(X,\Omega^2_X)$ and $P\in C$,
put $\omega(P)=\frac{\omega}{dz\wedge dv}(Q)$
for an arbitrary point $Q\in E_P$ and identify $\omega$ with the
pullback of a tensor on $C$. Denote by $\sigma\subset X$
the given zero-section of $f:X\to C$.

\begin{proposition}\label{inject}
  Suppose that $0\ne\omega\in H^0(X,\Omega^2_X)$
  and that $(\omega)_0$ is disjoint from $Z$. Then the
  cup product
  $\omega:H^1(X,\Theta_X(-\log\s))\to H^{1,1}_{prim}(X)$
  is injective.
  \begin{proof} Until now we have identified the tangent
    space $T_{\sje}(X)$ with a line bundle on $Z$; 
    we can also identify it with $H^1(X,\Theta_X(-\log\s))$.
    Since $\Theta_X(-\log\s)\otimes\Omega^2_X$ is
    isomorphic to $\Omega^1_X(\log\s)(-\s)$ we get a short exact
    sequence
    $$0\to \Theta_X(-\log\s)\to \Omega^1_X(\log\s)(-\s)
    \to \bigoplus_{P\in(\omega)_0}\sF_P\to 0,$$
    where $\sF_P$ is a rank 2 vector bundle on $E_P$
    that fits into a short exact sequence
    $$0\to\sO_{E_P}(-\s_P)\to\sF_P\to\sO_{E_P}\to 0.$$
    The coboundary map $H^0(E_P,\sO_{E_P})\to
    H^1(E_P,\sO_{E_P}(-\s_P))$ is identified with the Kodaira--Spencer map,
    so is an isomorphism from our assumption about $(\omega)_0$.
    So $H^0(E_P,\sF_P)=0$ and then the homomorphism
    $$H^1(X, \Theta_X(-\log\s))\to H^1(X, \Omega^1_X(\log\s)(-\s))$$
    is injective. This homomorphism factors through
    the cup product with which we are concerned,
    and the proposition is proved.
  \end{proof}
\end{proposition}

Write $b_{p,\u0}^{(j)}=b_p$. Then
$\omega^{(j)}=\sum b^{(j)}_{p\ge 0}z^pdz\wedge dv$
and $b^{(j)}_{0}=\omega^{(j)}(a)$.

\begin{theorem}\label{above}
  \part[i] $\eta_i^{(j)}\in H^0(X,\Omega^2_X(iE_a)) _{2^{nd}\ kind}$.
  \part[ii] $\eta_i^{(j)}=
  \sum_{p\ge 0}(i-p-1) b^{(j)}_pz^{p-i}dz\wedge dv$
  modulo $H^0(X,\Omega^2_X)$.
  \part[iii] The classes
  $[\eta_2^{(1)}],\ldots ,[\eta_2^{(h)}]$ span a
  line $\C[\eta_{a,2}]$ in $Fil^1/Fil^2$
  where $\eta_{a,2}\in
  H^0(X,\Omega^2_X(2E_a)) _{2^{nd}\ kind}$.
  \part[iv] More generally, for each $\ell$
  with $2\le \ell\le e(a)$, there exists a form
  $\eta_{a,\ell}\in H^0(X,\Omega^2_X(\ell E_a)) _{2^{nd}\ kind}$
  such that the classes
  $[\eta^{(j)}_k]_{j\in [1,h],k\in[2,\ell]}$
  span an $(\ell-1)$-dimensional subspace
  $\C.\{[\eta_{a,2},\ldots ,[\eta_{a,\ell}]\}$
  of $H^{1,1}_{prim}(X)$.
  
  \begin{proof} This follows from further consideration
    of the expansion (\ref{expand}). To begin, write
    $b^{(j)}_{p,\u0}=b^{(j)}_0.$ Then we get
    $$\omega^{(j)}(\ut)=
    \sum_{p,\uq}b^{(j)}_{p,\uq}z^p
    \left(1+\sum_i(i-p-1)t_iz^{-i}\right)dz\wedge dv,$$
    modulo $(\ut).H^0(X,\Omega^2_X)$,
    so that
    $$\omega^{(j)}=\sum_{p}b^{(j)}_{p}z^pdz\wedge dv\
    \textrm{and}\
    \eta^{(j)}_i=\sum_{p}
    (i-p-1)b^{(j)}_{p}z^{p-i}dz\wedge dv$$
    where the second equality holds modulo $H^0(X,\Omega^2_X)$.
    Observe that
    $$\dim H^0(X,\Omega^2_X(iE_a))\le h+i\ \textrm{and}\ 
    H^0(X,\Omega^2_X(E)) _{2^{nd}\ kind}=H^0(X,\Omega^2_X),$$
    so that
    $\dim H^0(X,\Omega^2_X(iE_a))_{2^{nd}\ kind}\le h+i-1.$
    Then inspection of these coefficients and the linear
    independence provided by Proposition \ref{inject}
    are enough to prove the theorem.
  \end{proof}
\end{theorem}

We can restate all this in more intrinsic terms, as follows.
Let $\uom$ denote the vector
$[\omega^{(1)},\ldots ,\omega^{(h)}]$
and $\uom_{(i)}$ its $i$'th derivative
with respect to $z$.

\begin{theorem}\label{deriv formula}
  For each $a\in Z$ and each $k=2,\ldots ,e(a)$
  there is an explicit meromorphic $2$-form
  $$\eta_{a,k}\in H^0(X,\Omega^2_X(kE_a))_{2^{nd}\ kind}$$
  such that
  $[\eta_{a,k}]$ lies in $Fil^1$
  and the image of
  $\bnabla_{\partial/\partial t_k}:Fil^2\to
  Fil^1$
  is spanned by the classes $[\eta_{a,2}],\ldots, [\eta_{a,k}]$.
  The kernel of $\bnabla_{\partial/\partial t_k}$
  contains $H^0(X,\Omega^2_X(-(k-1)E))$.

  Moreover, given the identification
  $H^0(X,\Omega^2_X)\cong H^0(X,\Omega^2_X)^\vee$
  provided by the basis $\uom$,
  $\bnabla_{\partial/\partial t_k}$
  is, as an element of
  $H^0(X,\Omega^2_X)\otimes Fil^1$,
  a linear combination of
  the tensors
  $\uom(a)\otimes[\eta_{a,k}],\ldots
  ,\uom_{(k-2)}(a)\otimes[\eta_{a,2}]$.

  Finally, If $h+q-1\ge e(a)$ then 
  $\bnabla_{\partial/\partial t_k}$
  is of rank $k-1$.
  \begin{proof}
    As remarked, this is, except
    for the final statement, nothing more than a
    restatement in intrinsic terms of what we have just proved.
    For the final part we need to know that the vectors
    $\uom(a),\ldots, \uom_{(e(a)-2)}$ are linearly
    independent. This follows from the cohomology
    of the exact sequence
    $$0\to\sO_C(K_C+L-(e(a)-1)a)\to
    \sO_C(K_C+L)\to
    \sO_{(e(a)-1)a}(K_C+L)\to 0$$
    and the assumption that $h+q-1\ge e(a)$.
  \end{proof}
\end{theorem}

\begin{theorem}\label{linearize}
  \part[i] If $a,b\in Z$ are distinct then, for all $v\in L_a$
  and for all $\eta_{b,k}$ with
  $k\in [2,e(b)]$,
  $\bnabla_v\eta_{b,k}$
  is a linear combination of the classes
  $\{[\eta_{a,i}]\}_{i\in [2, e(a)]}$ and
  $\{[\eta_{b,j}]\}_{j\in [2,e(b)]}$.
  In particular, $\bnabla_v\eta_{b,k}$
  lies in $Fil^1$.

  \part[ii] If $a,b\in Z$ are distinct then the
  classes $[\eta_{a,i}]$ and $[\eta_{b,k}]$ are orthogonal.

  \part[iii] The classes $[\eta_{a,i}]$ form a basis
  of $H^{1,1}_{prim}(X)$ as $a$ runs over the points of
  $Z$ and $i$ runs from $2$ to $e(a)$.
  
  \part[iv] If $Z$ is reduced then the classes $[\eta_{a,2}]$
  form an orthogonal
  basis of
  $H^{1,1}_{prim}(X)$ as $a$ runs over the points of $Z$.
  \begin{proof} \DHrefpart{i}:
    We use the same family $\sX\to\sC\to\Delta^{e(a)-1}$
    corresponding to the point $a$
    as before. This variation is constant outside a small \nbd
    of $a$ and the morphism to $\bsell$ is also
    constant outside this \nbd. So the ramification
    locus $Z$ is also constant there.

    The point $b$ moves in a family of points $b(\ut)\in\sC_{\ut}$,
    each $b(\ut)$ being of constant multiplicity $e(b)-1$
    in the ramification
    divisor $Z_{\ut}\subset \sC_{\ut}$. The fibre
    $E_b$ in $X$ moves
    in a family $\sE_{b(\ut)}$.
    Then, for each $k\in [2, e(b)]$,
    $$\eta_{b(\ut),k} =\eta_{b,k} + \sum t_i\nabla_{\partial/\partial
      t_i}\eta_{b(\ut),k}$$
    modulo $t^2$.
    Similarly to what we did before, we write
    $$\eta_{b(\ut),k} = \Res_{\sX_{\ut}}
    \left(
      \begin{array}{c}H_{b(\ut),k}\\
        {\ut}\\
      \end{array}
    \right)
    $$
    for some meromorphic $(e(a)+1)$-form
    $$H_{b(\ut),k}\in
    H^0(\sX, \Omega^{e(a)+1}_{\sX}(k\sE_{b(t)})).$$
    Then the same kind of calculation in terms of a
    power series expansion as before shows that
    $\nabla_{\partial/\partial t_i}\eta_{b(\ut),k}$
    is an element of
    $H^0(X,\Omega^2_X(iE_a + kE_b))$
    whose residues along both $E_a$ and $E_b$
    are zero.
    Therefore $\nabla_a(\eta_{b,k})$ is,
    modulo $H^{2,0}(X)$,
    a linear combination as described.

    \DHrefpart{ii}: Choose an element $\omega$
    of $H^{2,0}(X)$ that does not vanish along $E_a$.
    So $\nabla_{\partial/\partial t_i}\omega$ is a non-zero multiple
    of $\eta_{a,i}$. We can assume that
    $\nabla_{\partial/\partial t_i}\omega=\eta_{a,i}$.
    Now
    $\langle \omega, \eta_{b,k}\rangle =0$,
    since $H^{2,0}(X)$ is orthogonal to $H^{1,1}(X)$,
    so that
    $$0=
    \langle \nabla_{\partial/\partial t_i}\omega,\eta_{b,k}\rangle +
    \langle \omega, \nabla_{\partial/\partial t_i}\eta_{b,k}\rangle
    =\langle \eta_{a,i},\eta_{b,k}\rangle,
    $$
    since, as we have just proved,
    $\nabla_{\partial/\partial t_i}\eta_{b,k}\in Fil^1$.

    \DHrefpart{iii} and \DHrefpart{iv} follow from
    the linear independence of the $[\eta_{a,i}]$ and the
    fact that there are $N$ of them,
    where $N=\dim H^{1,1}_{prim}(X).$
  \end{proof}
\end{theorem}

\begin{remark} The fact that the $\eta_{a,i}$ form a basis
  of $H^{1,1}_{prim}(X)$ follows from Theorem 3.24
  and Remark 3.29 of
  \cite{CZ}, or from Proposition \ref{inject}.
  However, the orthogonality seems to be new.
\end{remark}

Until now $f:X\to C$ has been a point in $\sje^{gen}$.
Now, however, allow $X$ to have $A_1$-singularities,
so that $f:X\to C$ is defined by a classifying morphism
$\phi:C\to\bsell$ that is simply ramified over $j=\infty$.
Let $\tX\to X$
be the minimal model, so that $\tX$ has singular
fibres of types $I_1$ and $I_2$
and $X\to C$
is the relative canonical model of $\tX\to C$.
Suppose that
$D_1,\ldots, D_r$ are the exceptional
$(-2)$-curves on $\tX$.
Suppose that $a_1,\ldots, a_r$ are the points in $Z$
lying over $j=\infty$ and that $a_{r+1},\ldots, a_N$
are the other points of $Z$.
For $a=a_i$ with $i\le r$ define $[\eta_a]=[D_i]$.
The surface $\tX$ is a point in the stack $\tsje$
whose points are minimal models of points of
$\sje$. We can extend Theorem \ref{linearize}
to this context, as follows. For simplicity
we state it with the assumption that $Z$ is reduced.

\begin{theorem}\label{lin extra}
  Suppose that $Z$ is reduced. Put $\eta_{a,2}=\eta_a$
  for each $a\in Z$.
  
  \part[i]
  Each point $a$ in $Z$
  defines a line $v_a$ in the tangent space
  $H^1(\tX,T_{\tX})=T_{\sje}\tX$ such that
  the covariant derivative
  $$\bnabla_a:H^2(\tX,\Omega^2_{\tX})
  \to H^{1,1}_{prim}(\tX)$$
  is proportional to
  the rank one tensor $\omega_a^\vee\otimes \eta_a$.
  \part[ii] The classes $[\eta_a]_{a\in Z}$ form an orthogonal
  basis of $H^{1,1}_{prim}(\tX)$.
  \begin{proof} We only need to prove \DHrefpart{i} when
    $a=a_i$ for $i\le r$.
    Regard the surface $\tX$ as the specialization of a surface in
    $\sje^{gen}$. It follows by continuity that
    $\bnabla_a=\omega_a^\vee\otimes \theta_a$
    for some class $\theta_a$. To see that
    $\theta_a$ is proportional to $[D_a]$
    we argue as follows.

    Suppose that $\G=\sum \G_j$ is a configuration of
    $(-2)$-curves on a smooth surface $\tX$ that
    contracts to a single du Val singularity $P\in X$.
    Then we identify $H^2(\tX,\Omega^2_{\tX})
    =H^0(X,\omega_X)$ and put
    $$V(P)=H^0(X,\omega_X)/\mathfrak m_P H^0(X,\omega_X) .$$
    There is a short exact sequence
    $$0\to H^1(\tX,T_{\tX}(-\log \G))\to
    H^1(X,T_X)\to \oplus H^1(\G_j,\sN_{\G_j/\tX})$$
    which gives a commutative diagram
    $$\xymatrix{
      {H^1(\tX,T_{\tX})\otimes
        H^0(X,\omega_X)}\ar[r]\ar[d]
      &{H^1(\tX,\Omega^1_{\tX})}\ar[d]\\
      {H^1(\tX,T_{\tX})/H^1(\tX,T_{\tX}(-\log \G))\otimes V(P)}
      \ar[r] &{\oplus \C[\G_j]^\vee}.
    }$$
    Taking $\G=D_a$ leads to the fact
    that $\theta_a$ is proportional
    to $[D_a]$.

    \DHrefpart{ii} is proved as in Theorem \ref{linearize}.
  \end{proof}
\end{theorem}
\end{section}
\begin{section}{A local Schottky theorem}\label{section linearize}
    In this section we use the coincidence that
    $\dim H^{1,1}_{prim}=\dim \sje$ to put
    further structure on the derivative
    of the period map.
    
    The vector spaces $H^{1,1}_{prim}(X)$ fit together
    into a vector bundle $\sH=\sH^{1,1}_{prim}$ on $\sJE^{gen}$.
    We shall restrict attention to the open
    substack $\sJE^{ss}$ of $\sje^{gen}$ over which the universal
    ramification divisor $\sZ$ is {\'e}tale,
    so that, under the projection
    $\rho:\sZ^{ss}\to\sJE^{ss}$, the sheaf $\rho_*\sO_{\sZ^{ss}}$
    is a sheaf of semi-simple rings on $\sJE^{ss}$.

    \begin{proposition}\label{naturally}
      On $\sJE^{ss}$ the vector bundle
      $\sH$ is naturally
      a line bundle $\sB$ on $\sZ^{ss}$.
      \begin{proof}
        Essentially, we do this componentwise,
        using the orthogonal basis
        of $H^{1,1}_{prim}(X)$ that is
        provided by the classes $([\eta_a])_{a\in Z}$
        described in
        the previous section. 

        Let $x\in\sJE^{ss}$, and choose an analytic \nbd
        $U$ of $x$ such that $\rho^{-1}(U)$
        is a disjoint union
        $$\rho^{-1}(U)=\sqcup_{a\in Z} U_a$$
        of copies of $U$.

        Suppose $h\in \rho_*\sO_{\sZ}$. So
        $h\vert_{U_a}$ is identified with a holomorphic function
        $h_a$ defined on $U$. If
        $s\in \G(U,\sH)$, then we can write
        $s=\sum f_a[\eta_a]$
        where $f_a$ is a function on $U$.
        Now define $hs$ by the formula
        $$hs=\sum h_af_a[\eta_a].$$
      \end{proof}
    \end{proposition}

    \begin{proposition}\label{pairing} There is
      a perfect $\sO_{\sZ^{ss}}$-bilinear
      pairing $\b:\sB\times \sB\to \sO_{\sZ^{ss}}$
      such that $Tr(\rho_*\b)$ is the intersection product on $\sH$.
      \begin{proof} 
        We define the pairing $\b$ in terms
        of the notation used in the proof of Proposition
        \ref{naturally}, as follows:
        $$\b(s,t) =\b\left(\sum_a f_a[\eta_a],\sum_a g_a[\eta_a]\right)
        =\sum f_ag_a([\eta_a]\cup [\eta_a]).$$
        It is clear that $\b$ is perfect.
      \end{proof}
    \end{proposition}

    Recall that, tautologically, the tangent bundle
    $T_{\sJE^{gen}}$ is naturally a line bundle $\sS$ on $\sZ$.
    That is, $T_{\sJE^{gen}}=\rho_*\sS$. The derivative
    $per_*$ of the period map is an $\sO_{\sJE^{gen}}$-linear map
    $$per_*:\rho_*\sS \to \sH om_{\sJE^{gen}}(\sH^{2,0},\rho_*\sB).$$
    Define
    $$\sP=\Omega^1_{\sC/\sJE^{gen}}\otimes\Phi^*M,$$
    where $\sC$ is the pull back of the universal curve over $\sM_q$
    to $\sJE^{gen}$ and $\Phi:\sC\to \bsell$
    is the classifying morphism.
    Then
    $$\Omega^2_{\sX/\sJE^{gen}}\cong F^*\sP,$$
    where $F:\sX\to \sC$ is the universal Jacobian
    surface. There is an evaluation map
    $$eval: \rho^*\sH^{2,0}\to \sP\vert_{\sZ}$$
    and this map is surjective.

    \begin{proposition} On the locus $\sJE^{ss}$
      there is a factorization of $per_*$ given by a
      a commutative diagram
      $$\xymatrix{
        {\rho_*\sS} \ar[r]^>>>>>{\widetilde{per_*}}\ar[dr]_{per_*}&
      {\rho_*\sH om_{\sZ^{ss}}(\sP\vert_{\sZ^{ss}},\sB)}\ar[d]^>>>>>>{eval^\vee}\\
    &  {\rho_*\sH om_{\sZ^{ss}}(\rho^*\sH^{2,0},\sB)}
            }
    $$
    where $\widetilde{per_*}$ is $\rho_*\sO_{\sZ^{ss}}$-linear.
    \begin{proof}
      This is a consequence of part (3) of Theorem \ref{deriv formula}.
      \end{proof}
    \end{proposition}

    We have proved almost all of the next result.
    It is a local Schottky theorem in that it gives
    a precise description of the image of the tangent
    bundle to moduli under the period map.
  
    \begin{theorem}\label{local schottky} (Local Schottky)
      There is an isomorphism
    $${\widetilde{\widetilde{per_*}}}:\sS \to
    \sH om_{\sZ^{ss}}(\sP\vert_{\sZ^{ss}},\sB)$$
    of line bundles on $\sZ^{ss}$ such that
    $$\widetilde{per_*}=\rho_*{\widetilde{\widetilde{per_*}}}.$$
    \begin{proof} It is clear from the preceding discussion
      that there is a  homomorphism
      ${\widetilde{\widetilde{per_*}}}$
      of line bundles on $\sZ^{ss}$ with the desired properties.
      Since $per_*$ is injective on fibres,
      so is ${\widetilde{\widetilde{per_*}}}$, and we are done.
    \end{proof}
  \end{theorem}
  \end{section}
\begin{section}{Recovering $\Ram_\phi$ and $C$ from the infinitesimal
    period data}\label{recover}
  {\bf{For the rest of this paper we make the assumptions (\ref{(***)})}}.
  
  Suppose that $f_i:X_i\to C_i$ are two points of $\sJE^{ss}$ that have
  equal geometric genus $h$ and irregularity $q$, and that each
  classifying morphism $\phi_i:C_i\to\bsEll$ is simply ramified.
  Say $\Ram_{\phi_i}=Z_i$ and $\deg Z_i=N$, so that $N=10h+8(1-q)$.
  Recall that $L_=\phi_i^*M$.
  Our assumptions
  ensure that $Z_i$ is reduced and the linear system
  $\vert K_{X_i}\vert = f_i^*\vert K_{C_i}+L_i\vert$
  has no base points and embeds $C_i$ in $\P^{h-1}$
  as the canonical model of $X_i$. Note that
  $\deg C_i = h+q-1$.
  \begin{theorem}\label{7.1} Assume that both surfaces have the same
    IVHS. Then $C_1=C_2$ and $Z_1=Z_2$.
    \begin{proof} Our assumptions mean that the two surfaces give the
      same point and
      the same tangent space under the period map after each
      cohomology
      group $H^2(X_i,\Z)$ has been appropriately normalized. In
      particular,
      such a normalization gives a normalized basis
      $\underline{\omega}_i=[\omega^{(1)}_i,...,\omega^{(h)}_i]$
      of $H^0(X_i,\Omega^2_{X_i})$
      for each $i$. We can therefore regard the curves $C_i$
      as embedded in the same projective space $\P^{h-1}$.
      That is, any point $P$ of $C_i$ is identified with the point
      $\underline{\omega}_i(\tP)$ in $\P^{h-1}$, where $\tP\in X$
      is any point of $X$ that lies over $P$. The basis
      $\underline{\omega}_i$
      also gives an identification of
      $H^0(X_i,\Omega^2_{X_i})$ with its dual.

      We show first that $Z_1=Z_2$. We then
      deduce,
      by quadratic interpolation,
      that $C_1=C_2$.

      By the results of Section \ref{schiffer},
  especially Theorem \ref{deriv formula},
  the tangent space $T_{\sJE}X_1$ is of dimension
  $N$
  and its image in
  $H^0(X_1,\Omega^2_{X_1})^\vee\otimes H^{1,1}_{prim}(X_1)$
  is a direct sum $\bigoplus_{a\in Z_1} L_{a}$ where
  $L_{a}$ is a line described by Theorem \ref{deriv formula}.
  We make the identifications
  $$H^0(X_1,\Omega^2_{X_1})=
  H^0(X_2,\Omega^2_{X_2})=U\
  \textrm{and}\ H^{1,1}_{prim}(X_1)=H^{1,1}_{prim}(X_2)=V$$
  and assume that the tangent spaces
  $T_{\sJE}X_1$ and $T_{\sJE}X_2=\bigoplus_{b\in Z_2}M_{b}$
  are equal as subspaces of $U\otimes V$.
  We know that $L_{a}$ is spanned by the rank $1$ tensor
  $u^{(i)}_0\otimes v^{(i)}_0$ (the image
  of $\bnabla_{\partial/\partial t_2}$).
  Similarly $M_{b}$ is spanned by a rank 1 tensor
  $p^{(j)}_0\otimes q^{(j)}_0$.

  \begin{lemma}\label{unscrew}\label{bang}
    Suppose that
    \begin{enumerate}
    \item $U$ and $V$ are vector spaces such that
      $\dim U=h$ and $\dim V=N$,
    \item $u_1,...,u_N\in U$,
    \item no two of the $u_i$ are linearly dependent,
    \item $v_1,...,v_N\in V$ and form a basis of $V$,
    \item $\xi\in U\otimes V$ and there is a linear relation
      $$\xi=\sum_1^N\lambda_k u_k\otimes v_k$$
      and, finally,
    \item $\xi=x\otimes y$, a tensor
      of rank $1$.
    \end{enumerate}
    Then there is a unique index $i$
    such that $\xi$ is proportional to $u_i\otimes v_i$.
    \begin{proof}
      Let $(v_i^\vee)$ be the dual basis of $V^\vee$.
      There is an index $j$ such that
      $\langle y, v_j^\vee\rangle\ne 0$.
      Then $$\langle y, v_j^\vee\rangle x=
      \sum \lambda_k\delta_{jk}u_k=\lambda_ju_j,$$
      so $x$ is a  multiple of $u_j$.
      Since no two of the $u_i$ are linearly dependent,
      this index $j$ is unique, so $y=\mu_jv_j$
      and then $x\otimes v_j=\sum \lambda'_ku_k\otimes v_k.$
    \end{proof}
  \end{lemma}

      Say $Z_1=\{a_1,...,a_N\}$ and
      $Z_2=\{a'_1,...,a'_N\}$. Take $U= H^{2,0}(X_1)$
      and $V=H^{1,1}_{prim}(X_1)$ and write
      \begin{eqnarray}\label{notation}
      {x_k={\underline{\omega}}_1(a_k),}\ 
      {x'_k={\underline{\omega}}_2(a'_k),}\ 
      {y_k=[\eta_{1,a_k}]}\ {\textrm{and}}\ 
                                            {y'_k=[\eta_{2,a'_k}].}
      \end{eqnarray}
            From Theorem \ref{deriv formula} and
      the assumption that the IVHS
      of the two surfaces $X_1$ and $X_2$ are isomorphic,
      the tensors $x'_l\otimes y'_l$ span the same $N$-dimensional
      subspace
      of $U\otimes V$ as do the tensors $x_k\otimes y_k$.
      In particular, each $x'_l\otimes y'_l$ is a linear
      combination of the tensors $x_k\otimes y_k$.
      Then, by Lemma \ref{unscrew}, for each index $l$ there
      is a unique index $m$ such that
      $x_l$ is proportional to $x'_m$.
      That is, $a_l=a'_m$, and therefore $Z_1= Z_2=Z$,
      say.

      By the assumptions (\ref{(***)})
      each $C_i$ is non-degenerately embedded in
      $\P^{h-1}$ by a complete linear system, and
      $\deg C_i\ge 2q+2$, so that,
      by the results of \cite{Mu} and \cite{SD},
      $C_i$ is an intersection of quadrics. Since
      $$\deg Z >2\deg C_i,$$
      again by (\ref{(***)}),
      each $C_i$ equals the intersection of the quadrics through
      $Z$, so that $C_1=C_2$.
          \end{proof}
        \end{theorem}

        We can reformulate this as follows.

        \begin{theorem}\label{reform} If the ramification
          divisor $\Ram_\phi$ of the classifying morphism
          $\phi:C\to\bsell$
          is reduced
          then the IVHS of the surface $X$
          determines the base curve $C$, the
          divisor $\Ram_\phi$
          and the line bundle $L=\phi^*M$.
          \begin{proof} What remains to be done
            is to show that we can recover the bundle $L$.
            For this, observe that our argument has shown that the IVHS
            of $X$ determines the embedding $i:C \inj \P^{h-1}$
            as the canonical model of $X$ and that
            $i^*\sO(1)=\sO_C(K_C+L).$
          \end{proof}
        \end{theorem}
          
        \begin{remark}
          \part[i] 
    Theorems \ref{7.1} and \ref{reform}
    lead to the problem of trying to
    recover a morphism $\phi:C\to\bsell$ from knowledge of $C$,
    the divisor $\Ram_\phi$ and the line bundle $L$.
    However, it is impossible to do more than
    recover $\phi$ modulo the action of the
    automorphism group $\GG_m$ of $\bsell$.
    In this direction we shall prove Proposition \ref{8.9}.
    \part[ii] If instead of $\Ram_\phi$ being reduced
    we assume only that $h+q-1\ge e(a)+e(b)-1$
    for all $a,b\in\Ram_\phi$ then a refinement of the
    argument given here shows that Theorem \ref{reform}
    still holds.
      \end{remark}
    \end{section}
    \begin{section}{The structure of the tangent bundle to $\sse^{gen}$}
      Recall from Section \ref{comparison}
      that we have a morphism $\pi:\bsm_1\to\bsell$
      that gives, locally on $\bsell$, an isomorphism
      $\bsm_1\to B\sG$, so that the tangent complex
      $T^\bullet_{\bsm_1}$
      is locally isomorphic to the $2$-term complex
      $$0\to \pi^*M^\vee[1] \to \pi^*T_{\bsell}[0]\to 0$$
      whose differential is zero.
      The morphism $\pi$ determines a morphism
      $$\Pi:\sse^{gen} \to\sJE^{gen}.$$
      Suppose that $f:X\to C$ corresponds to
      $\psi:C\to\bsm_1$ and maps under $\Pi$
      to $Y\to C$. Say $\phi=\pi\circ\psi$.
      Observe that $\Ram_\phi=\Ram_\psi=Z,$
      say.
      The description just given
      of $T^\bullet_{\bsm_1}$ leads to a short
      exact sequence
      $$0 \to H^1(C, -L) \to T_{\sse}X \to T_{\sje}Y \to 0,$$
      where $L=\phi^*M$. The subobject
      $H^1(C,-L)$ in this sequence
      is identified with the tangent space $T_{\Pi^{-1}(Y)}X$
      to the fibre of $\Pi$ though $X$.
      So
      $$T_{\Pi^{-1}(Y)}(X)\cong H^1(C,-L)\cong H^0(C,K_C+L)^\vee,$$
      which is in turn naturally isomorphic to both
      $H^{2,0}(X)^\vee$ and to $H^{2,0}(Y)^\vee$.
      Let $\xi$ denote the class of a fibre of $X\to C$;
      then
      the period map gives a commutative diagram
      with exact rows
      $$\xymatrix{
        {0} \ar[r]& {T_{\Pi^{-1}(Y)}(X)}\ar[d]\ar[r] &
        {T_{\sse}(X)}\ar[d]_{per_{X,*}}\ar[r] &
        {T_{\sje}(Y)}\ar[r]\ar[d]^{per_{Y,*}}&
        {0}\\
        {0}\ar[r]& {H^{2,0}(X)^\vee\otimes\xi}\ar[r]&
        {H^{2,0}(X)^\vee\otimes\xi^\perp}\ar[r]&
        {H^{2,0}(Y)^\vee\otimes H^{1,1}_{prim}(Y)}\ar[r]&
          {0.}
        }
        $$

        \begin{theorem}
          If $Z$ is reduced and disjoint from the locus
          $j=\infty$ then the IVHS of the surface $X$
          determines the ramification locus $Z$, the
          base curve $C$ and the line bundle $\phi^*M$.
          \begin{proof}
            Put $H^{2,0}(X)^\vee=H^{2,0}(Y)^\vee=U$
            and $\xi^\perp=V$, so that
            $\dim U=h$ and $\dim V=N+1$.
            The image of $T_{\sse}(X)$ under $per_{X,*}$ is an
            $(N+h)$-dimensional subspace
            $W$ of $U\otimes V$ such that $W$
            contains a subspace $U\otimes\xi$ and
            $W/U\otimes\xi$,
            which is the image of
            $T_{\sje}(Y)$ under $per_{Y,*}$,
            is spanned by rank one tensors
            $$x_1\otimes y_1,\ldots,x_N\otimes y_N$$
            where the vectors $y_1,\ldots,y_N$
            form a basis of $V/\xi =H^{1,1}_{prim}(Y)$.
            Let $\pi:V\to V/\xi$ denote the projection.

            Suppose that the same IVHS arises also from
            another surface $X'$. Then there is a vector
            $\xi'\in V$ that arises from $X$
            such that $(1_U\otimes\pi)(U\otimes\xi')$
            is a subspace of $U\otimes V/\xi$
            and lies in the subspace of $U\otimes V/\xi$
            that is spanned by the $x_i\otimes y_i$.
            However, by Lemma \ref{bang},
            the $x_i\otimes y_i$ are the only rank one tensors
            in $U\otimes V/\xi$, so that
            $(1_U\otimes\pi)(U\otimes\xi')=0$ and
            therefore $\xi'$ is proportional to $\xi$.
            Consider the vectors $x_i',y_i'$ that arise
            from $X'$; then the tensors $x_i\otimes y_i$
            and $x_i'\otimes y_i'$ lie in the same
            vector space $U\otimes V/\xi$, and
            then we can use Lemma \ref{bang}
            again to conclude the proof.
          \end{proof}
        \end{theorem}
        
    Compare the case of $\sM_q$:  over
    the non-hyperelliptic locus Schiffer variations
    give a cone structure in the tangent bundle, where
    at each point $C$ of $\sM_q$
    the corresponding cone is the cone over the bicanonical model
    of $C$, the generators of the cone map,
    under the period map,
    to tensors (quadratic forms) of rank $1$
    and, again, these account for all the rank $1$
    tensors in the image.
  \end{section}
  \begin{section}{Recovering information
      from $C$ and $Z$}
    \label{towards}
    We shall show
    (Proposition \ref{8.9}) that
    $\phi_1:C\to\bsell$ is generic and
    if $\phi_2:C\to\bsell$ is another morphism
    such that $\Ram_{\phi_1}=\Ram_{\phi_2}$
    and $\phi_1^*M\cong \phi_2^*M$, then
    $\phi_1$ and $\phi_2$ are equivalent modulo
    the action of $\Aut_{\bsell}=\GG_m$ provided that
    also there are sufficiently many points
    $a_i\in Z$ such that $\phi_1(a_i)$
    is isomorphic to $\phi_2(a_i)$. Therefore
    an effective form of 
    generic Torelli holds for Jacobian elliptic surfaces
    modulo this action of $\GG_m$.

    To begin, we rewrite some results of Tannenbaum \cite{T}
    in the context of Deligne--Mumford stacks.
    Assume that $\sS$ is a smooth
    $2$-dimensional Deligne--Mumford stack (the relevant
    example will be $\sS=\bsell\times\bsell$),
    that $C$ is a smooth projective curve and that
    $\pi:C\to \sS$ is a morphism that factors as
    $$
        \xymatrix
        {
          {C}\ar[d]_{\pi'}\ar[dr]^{\pi}&\\
          {D}\ar[r]^>>>>>i&{\sS}
        }
    $$
    where $D$ is a projective curve with only cusps,
    $\pi'$ is birational and $i$ induces surjections
    of henselian local rings at all points. (We shall say that
    ``$\pi$ is birational onto its image''.) Then there is a
    conormal sheaf $\sN^\vee_D$, a line bundle on $D$,
    which is generated by the pull back under $i$ of the
    kernel $\sI_D$ of $\sO_{\sS}\to i_*\sO_D$. It
    fits into a short exact sequence
    $$0\to \sN^\vee_D\to i^*\Omega^1_{\sS}\to\Omega^1_D\to 0.$$
    So the adjunction formula (or duality for the morphism
    $i:D\to \sS$) gives an isomorphism
    $$\sN_D\buildrel\cong\over\to
    \omega_D\otimes i^*\omega_{\sS}^\vee.$$
    There is also a homomorphism
    $\sN_D\to T^1(D,\sO_D)$; let $\sN'_D$ denote its kernel.
    The space $H^0(D,\sN'_D)$ is the tangent space
    to the functor that classifies deformations
    of the morphism $i:D\to \sS$ that are locally trivial
    in the {\'e}tale (or analytic) topology.
    
        These conditions imply that there is a reduced effective
        divisor $R$ on $C$ such that
        $T_C(R)=T_{C}\otimes\sO_C(R)$ is the saturation
        of $T_{C}$ in $\pi^*T_{\sS}$. Define $\sN'_\pi$ by
        $$\sN'_\pi=\coker(T_C(R)\to\pi^*T_{\sS});$$
        this is a
        line bundle on $C$.

        Define $\sJ\subset \sO_{\sS}$ to be the Jacobian ideal
        of the ideal $\sI_D$. The next result is a slight variant
        of Lemma 1.5 of \cite{T}.

        \begin{lemma}\label{tann} $\pi'_*\sN'_\pi=\sN_D'$.
          \begin{proof} We have $\Omega^1_D=\sJ.\omega_D$,
            by direct calculation,
            and, from the definition of $\sN'_D$, we have
            $\sN'_D=\sJ.\sN_D$.

            From the definition,
            $\omega^\vee_{C}(R)\otimes
            \sN'_\pi\cong\pi^*\omega_{\sS}^\vee$.
            Let $\sC$ denote the conductor ideal. Then
            $\omega_{C}\cong\sC\otimes\pi'^*\omega_D$,
            so that
            $$\sN'_\pi\cong
            \pi'^*\omega_D\otimes\sC\otimes\sO_{C}(-R-K_{\sS}).$$
            Now $\sC=\sO_{C}(-R)$, by the nature of a cusp,
            and $\sJ.\sO_{C}=\sO_{C}(-2R)$ for the same reason.
            So
            $$\sN'_\pi\cong\pi'^*(\omega_D\otimes
            i^*\omega_S^\vee)\otimes\sO_{C}(-2R)$$
            and therefore
            $$\pi'_*\sN'_\pi\cong
            \sJ.(\omega_D\otimes\omega_{\sS}^\vee)\cong\sJ.\sN_D=\sN'_D.$$
          \end{proof}
        \end{lemma}

        \begin{corollary} $H^0(C,\sN'_\pi)$ is isomorphic
          to the tangent space of the deformation functor
          that classifies those deformations of the morphism
          $\pi:C\to \sS$ where the length of the $\sO_{C}$-module
          $\coker(\pi^*\Omega^1_{\sS}\to\Omega^1_{C})$
          is preserved.
          \begin{proof}
            By the lemma, $H^0(C,\sN'_\pi)$ is the tangent
            space to locally trivial deformations of the morphism
            $i:D\to \sS$. Since $D$ has only cusps, a deformation
            of $i$ is locally trivial if and only if it preserves the
            length of $\coker(\pi^*\Omega^1_{\sS}\to\Omega^1_{C})$.
          \end{proof}
        \end{corollary}

        \begin{proposition}\label{8.3} Fix a generic
          morphism
          $\phi_1:C\to\bsell$ and consider the
          morphisms $\phi_2:C\to\bsell$ such that
          $(\phi_1,\phi_2):C\to\bsell\times\bsell$ is birational onto
          its image.

          \part[i] If $q\ge 2$ there are no such morphisms $\phi_2$.
          \part[ii] If $q=1$ there are only finitely many such
          $\phi_2$.
          \part[iii] If $q=0$ then these morphisms $\phi_2$
          form a family of dimension at most $3$.
          \begin{proof} Assume that $\phi_2$ exists. Take
            $\bsell\times\bsell=\sS$ and $(\phi_1,\phi_2)=\pi$.
            Then there is a factorization
            $$
        \xymatrix
        {
          {C}\ar[d]_{\pi'}\ar[dr]^{\pi}&\\
          {D}\ar[r]^i&{\sS}
        }
        $$
        as before, and the divisor $R$ is
        $R=\Ram_{\phi_1}=\Ram_{\phi_2}$.
        We get
            \begin{eqnarray*}
              \deg \sN'_\pi&=&2\deg\phi_1^*T_{\bsell}+2q-2-N\\
                         &=& 2.24(h+1-q).\frac{5}{12}+2q-2-N\\
              &=&10h+10(1-q),
            \end{eqnarray*}
            \noindent so that
            $$h^0(C,\sN'_\pi)=10h+11(1-q)=N-(3q-3).$$
            But $\dim\sje=N,$
            and the proposition is proved.
          \end{proof}
        \end{proposition}

        \begin{corollary} Fix $r\ge 0$ and
          suppose that $\phi_1:C\to\bsell$
          is generic. Consider the morphisms
          $\phi_2:C\to\bsell$ such that
          $\Ram_{\phi_1}=\Ram_{\phi_2}=Z$,
          say,
          and there exist distinct points $a_1,\ldots,a_r\in Z$
          such that $\phi_1(a_i)$ is isomorphic to
          $\phi_2(a_i)$ for all $i=1,\ldots, r$
          and moreover
          $(\phi_1,\phi_2):C\to \bsell\times\bsell$
          is birational onto its image.

          Then no such $\phi_2$ exist provided that one of the
          following is true:
          $q\ge 2$; $q=1$ and $r\ge 1$; $q=0$ and $r\ge 4$.
          \begin{proof} The further constraints on $\phi_2$
            imply that, in the notation of the proof of Lemma
            \ref{8.3}, the curve $D$ has $r$ cusps that lie
            on the diagonal of $\sS$. Therefore there are
            $r$ further constraints on $\phi_2$ and the
            corollary follows.
          \end{proof}
        \end{corollary}

        \begin{lemma}\label{8.4} If $\phi:C\to\bsell$ is generic then there
          is no non-trivial factorization through a curve of the composite morphism
          $\g:C\to\P^1_j$.
          \begin{proof} Suppose that
            $$C\buildrel\beta\over\to\G\buildrel\a\over\to\P^1_j$$
            is a non-trivial factorization of $\g=\pi\circ\phi$.
            Say $\deg\a =a$ and $\deg\beta= b$.
            By the assumption of genericity, the divisor
            $Z=\Ram_\phi$is reduced and its image $B=\g_*Z$ in
            $\P^1_j$ consists of distinct points, none
            of which equals $j_4$ or $j_6$.

            Suppose $\a$ is branched over $y\in\P^1_j$
            and $y\ne j_4,j_6$. Then $\g^{-1}(y)\subset Z$,
            so that $Z\to B$ is not 1-to-1. So
            $\a$ is branched only over $j_4,j_6$,
            and so is the cyclic cover of $\P^1_j$
            that is
            branched at these two points and
            is of degree $a$.
            Since $\phi$ is {\'e}tale over $j_4,$
            and over $j_6$
            it follows that $a\vert 2$ and $a\vert 3$,
            contradiction.
          \end{proof}
        \end{lemma}
        
        Suppose that $\Delta$ is any irreducible curve
        of bidegree $(1,1)$ in $\P^1_j\times\P^1_j,$
        let $\tDelta'$ denote 
        the fibre product
        $$\tDelta'=\Delta\times_{\P^1_j\times\P^1_j}
        \bsell\times\bsell$$
        and let $\tDelta$ be the normalization
        of $\tDelta'$. Note that, generically,
        $\tDelta$ is isomorphic to $\Delta\times B(\Z/2\times\Z/2)$,
        so that $\deg(\tDelta\to\Delta)=\frac{1}{4}$.

        \begin{lemma}\label{8.8}
          $\deg T_{\tDelta}=5/24$
          if $\Delta$ contains the points
          $(j_4,j_4)$ and $(j_6,j_6)$
          and $\deg T_{\tDelta}\le 1/12$
          otherwise.
          \begin{proof} We consider separately the
            cases according to which of the special points
            $(j_4,j_4)$, $(j_4,j_6)$,
            $(j_6,j_4)$ and $(j_6,j_6)$
            lie on $\Delta$
            and summarize the results in
            Table \ref{Tab:table 2}. Each case depends on local
            calculation.
            For example, suppose that $(j_4,j_4)$ and
            $(j_6,j_6)\in\Delta$.
            Then $\tDelta'$ has a double point with an action of
            $(\Z/4)^2$ lying over $(j_4,j_4)$
            and a triple point with an action of $(\Z/6)^2$
            over $(j_6,j_6)$. In each case the branches are permuted
            transitively
            by the group, so that $\tDelta$ has one point with an
            action of $\Z/4\times\Z/2$ and another with an action
            of $\Z/6\times\Z/2$. So
            \begin{eqnarray*}
              \deg T_{\tDelta}&=&\deg T_{\Delta}.\deg B((\Z/2)^2)+
                                    (-\deg B((\Z/2)^2)+\deg B(\Z/4\times\Z/2))\\
              &+&
                  (-\deg B((\Z/2)^2)+\deg B(\Z/6\times\Z/2))=1/8+1/12=5/24.
                  \end{eqnarray*}
            \begin{table}[h]
              \begin{center}
              \captionof{table}{{\label{Tab:table 2}}}
              \begin{tabular}{|c|c|}
                \hline
                {Special points on $\Delta$}&{$\deg T_{\tDelta}$}\\
                               \hline
                {None}&{$-1/12$}\\
                {$(j_4,j_4)$}&{$1/24$}\\
                {One of $(j_4,j_6)$ and $(j_6,j_4)$}&{$0$}\\
                                {$(j_6,j_6)$}&{$1/24$}\\
                {$(j_4,j_6)$ and $(j_6,j_4)$}&{$1/12$}\\
                {$(j_4,j_4)$ and $(j_6,j_6)$}&{$5/24$}\\[1ex]
                \hline
              \end{tabular}
            \end{center}
          \end{table}
        \end{proof}
      \end{lemma}
      \begin{proposition}\label{8.9} 
        Suppose that $\phi_1:C\to\bsell$ is generic
        and that $\phi_2:C\to\bsell$ is another morphism
        such that $\Ram_{\phi_1}=\Ram_{\phi_2}=Z$
        and there are $r$ points $a_i\in Z$ such that
        $\phi_1(a_i)$ is isomorphic to $\phi_2(a_i)$
        for all $i$. Assume also that one of the following is true:
        $q\ge 2$; $q=1$ and $r\ge 1$; $q=0$ and $r\ge 4$.
        
        \part[i] 
        The image of the composite morphism
        $C\to\P^1_j\times\P^1_j$
        is a curve $\Delta$ of bidegree $(1,1)$.

        \part[ii] $\Delta$
        passes through the points
        $(j_4,j_4)$ and $(j_6,j_6)$.

        Suppose also that $\phi_1^*M\cong \phi_2^*M$.

        \part[iii] The
        morphism $(\phi_1,\phi_2): C\to\bsell\times\bsell$
        factors through
        the graph of an automorphism of $\bsell$.

        \part[iv] $\phi_1$ and $\phi_2$ are equivalent
        under the action of $\Aut_{\bsell}$.
        \begin{proof}
          \DHrefpart{i} follows from Proposition
          \ref{8.3} and Lemma \ref{8.4}.

          \DHrefpart{ii}: Suppose  this is false.
          Then, by Lemma \ref{8.8},
          $$\dim_{(\phi_1,\phi_2)}\Mor(C,\tDelta)
          \le 2\deg\phi_1\times\frac{1}{12}+1-q.$$
          But $\dim_{(\phi_1,\phi_2)}\Mor(C,\tDelta)
          \ge \dim_{\phi_1}\Mor(C,\bsell)$,
          contradiction.

          For \DHrefpart{iii} and \DHrefpart{iv}
          we can now use
          the $\GG_m$-action to ensure 
          that $\Delta$ is the diagonal. Then the two
          surfaces $X_i\to C$ have the same $j$-invariant,
          so that one is a quadratic twist of the other
          via a quadratic cover $\tC\to C$. Since
          both surfaces are semistable, this quadratic
          cover is {\'e}tale. Suppose it corresponds
            to the $2$-torsion class $P$ on $C$.
            From the identification of
            $\phi_i^*M$ with the conormal
            bundle of the zero section, it follows
            that $\phi_2^*M= \phi_1^*M + P$,
            so that $P=0$, and the result is proved.
           \end{proof}
         \end{proposition}
       \end{section}
       \begin{section}{Generic Torelli for Jacobian
           surfaces}\label{gen tor}
         Define $r(q)= 2$ if $q\ge 1$ and $r(0)=4$. Let $\tsje_{C,L}$
         (or $\sje_{C,L}$)
         denote the closed substack of $\tsje$ (or $\sje$)
         defined by the properties
         that the base of the elliptic fibration is $C$ and $\phi^*M$
         is isomorphic to $L$. Let $\tsje_r$ denote the substack of
        $\tsje$ with $r$ fibres of type $I_2$ and put
        $\tsje_{C,L,r}=\tsje_{C,L}\cap\sje_r$.

        Let $\sZ'\to\sje_r$ be the restriction of the universal
        ramification divisor $\rho:\sZ\to\sje$. Then $\sZ'$ has a
        closed substack $\sZ_r$ obtained by discarding the part
        of $\sZ$ that lies over $j=\infty$. 

        \begin{lemma}
          $\tsje_r$, $\sje_r$, $\tsje_{C,L,r}$, $\sje_{C,L,r}$ and $\sZ_r$
          are irreducible when $r=r(q)$.
          \begin{proof} We show first
            that $\sje_{C,L,r}$ is irreducible. For this,
            refer to the Weierstrass equation
            $$y^2=4x^3-g_4x-g_6$$
            of a Jacobian surface $X$ over $C$.
            Fix distinct points $P_1,\ldots, P_r\in C$.
            Then $X$ has fibres of type $I_2$ (or worse) over each
            $P_i$ if and only if 
            $$g_4^3-27 g_6^2= 2g_4g_6'-3g_6g_4'=0,$$
            where the prime denotes the derivative
            with respect to a local co-ordinate $z$ at $P_i$.
            (For this, we use the description of $\sZ$
            in terms of transvectants that is given
            in Section \ref{schottky} below.)
            Expand $g_n$ in terms of $z$ to get
            $$g_4=a_0+a_1 z+h.o.t.,\ g_6=b_0+b_1 z+ h.o.t.;$$
            then the condition on the fibres of $X$
            is expressed by a pair of equations
            (one pair for each $P_i$)
            $$a_0^3-27b_0^2=2a_0b_1 - 3a_1b_0=0.$$
            Since $r\le 4$ these equations define
            an irreducible subvariety $V_r$ of $H^0(C,4L)\times H^0(C,6L)$
            and the irreducibility of $\sje_{C,L,r}$ is proved. The
            irreducibility of $\sje_r$ follows at once, and then
            so too does the irreducibility of $\tsje_{C,L,r}$
            and $\tsje_r$.

            Since $\sZ_r\cap\sje_{C,L,r}$ is dominated by $V_r$,
            the irreducibility of $\sZ_r$ is also proved.
          \end{proof}
        \end{lemma}
\begin{theorem} 
\part[i] Suppose that $\phi_i:C\to\bsell$
is a morphism corresponding to the Jacobian surface
$f_i:X_i\to C$, for $i=1,2$. 
Assume that $\Ram_{\phi_1}=\Ram_{\phi_2}=Z$, say,
that $Z$ is reduced and that $\phi_1^*M=\phi_2^*M$. 
Suppose that there are $r=r(q)$ points $a_1,\ldots, a_r\in Z$
such that $\phi_i(a_k)$ lies over $j=\infty$
for $i=1,2$ and $k=1,\ldots, r$. Then $X_1$
is isomorphic to $X_2$.

\part[ii] A generic Jacobian elliptic surface
with $r(q)$ singular fibres of type $I_2$ is determined
by its IVHS.

\part[iii] Generic Torelli holds for Jacobian
surfaces with $r(q)$ singular fibres of type $I_2$.
\begin{proof} \DHrefpart{i} is a consequence
  of Proposition \ref{8.9}. \DHrefpart{ii} is
merely a restatement of \DHrefpart{i},
and then \DHrefpart{iii} follows in the usual way.
\end{proof}
\end{theorem}

Next, we use the MMP to prove a version
of the good reduction result from \cite{C1}
and \cite{C2}.

\begin{theorem} Suppose that $\sX\to\Delta$ is a
  semi--stable
  $1$-parameter degeneration of Jacobian elliptic surfaces
  over a fixed curve $C$ of genus $q$
  and that there is no monodromy on the second
  Betti cohomology of 
  $\sX_{\bar\eta}$. 
  Assume also that, under the period map,
  the image of $0\in\Delta$ is the primitive weight 2 Hodge structure of a
  Jacobian elliptic surface $g:V\to C$ such that $p_g(V)=p_g(\sX_{\bar\eta})$
  and this Hodge structure
  is irreducible.

  Then $\sX\to\Delta$ has good reduction.
  \begin{proof} We can follow the proof of Theorem 9.1 of \cite{SB}
    as far as the end of Lemma 9.9 of {\it{loc. cit}}. This leads to a model
    $\sX\to\Delta$ such that
    \begin{enumerate}
    \item $\sX$ has $\Q$-factorial canonical singularities
      and $\sX_0$ has slc singularities,
    \item there is a surface $S$ with a proper semi--stable morphism
      $g:S\to\Delta$ and a projective morphism $f:\sX\to S$ with only
      $1$-dimensional fibres such that $K_{\sX/\Delta}$ is the pullback
      under $f$ of a $g$-ample $\Q$-line bundle $L$ on $S$.
    \end{enumerate}
    Then $S_0$ is a tree, and $S_0=\sum C_i$ where $C_0\cong C$
    and $C_j\cong\P^1$ for all $j>0$. Since $S$ has only
    singularities of type $A$, each curve $C_i$ is a $\Q$-Cartier
    divisor on $S$. Let $f_i:X_i=f^{-1}(C_i)\to C_i$
    be the morphism induced by $f$. So
    $X_i$ is $\Q$-Cartier on $\sX$ and
    $K_{\sX/\Delta}\vert_{X_i}\sim f_i^*L_i$ for some
    $\Q$-line bundle $L_i$ on $C_i$
    of degree $\a_i>0$, where $\sum\a_i=h+q-1$.
    
    Since the Hodge structure on $H_{prim}^2(V)$ is irreducible
    and,
    from the absence of monodromy,
    $\sum p_g(X_i)\le p_g(V)$, there is a unique index
    $i_0$ such that the Hodge structure on $H^2(X_{i_0},\Q)$
    maps onto that on $H^2(V,\Q)$. It follows that
    $i_0=0$.

    Suppose that $X_0$ meets $X_1,\ldots, X_s$.
    Since $X_0$ is $\Q$-Gorenstein, and even Gorenstein
    in codimension one, the adjunction formula
    shows that there exist $\delta_1,\ldots, \delta_s\in\Q_{>0}$
    such that $K_{X_0}\sim_\Q f_0^*(\a_0-\sum\delta_i)$.
    
    Let $\g:\tX\to X_0$ denote the minimal
    resolution.
    Then there is an effective $\Q$-divisor $D$ on $\tX$
    such that
    $$K_{\tX}\sim_\Q
    \g^*K_{X_0}-D\sim_\Q
    (f_0\circ\g)^*(\a_0-\sum\delta_i)-D.$$
    So
    $$p_g(\tX)\le\lfloor \a_0-\sum\delta_i\rfloor +1-q
    \le \sum\a_j+1-q$$
    On the other hand,
    $$p_g(\tX)=p_g(V)=h=\sum\a_j+1-q.$$
    Therefore $s=0$, so that $\sX_0=X_0$,
    and $D=0$. Therefore $X_0$
    has only du Val singularities and
    the result is reduced to the existence of a
    simultaneous resolution for such singularities.
  \end{proof}
\end{theorem}

\begin{remark} There are degenerating families of elliptic surfaces
  $\sX\to\sC\to\Delta$ with no monodromy on $H^2$ but where
  $\sC$ has bad reduction. For example, take an elliptic
  K3 surface $X\to\P^1$ with two isomorphic fibres of type $\bD_4$.
  Then this can be plumbed to itself to give such a family:
  the generic fibre has $p_g=q=1$ and $\sC\to\Delta$
  is the Tate curve.
\end{remark}

\begin{corollary}\label{corollary} The period map
  $[per]:[\sje]\to [D/\Gamma]$
  of geometric quotients is proper and one-to-one
  over the image of the generic point of the
  locus of surfaces with $r(q)$
  fibres of type $I_2$.
  \noproof
\end{corollary}

\begin{lemma} Suppose that $V\inj Y$ is a closed
  embedding of normal analytic spaces. Suppose also
  that $H$ is a finite group and $G$ a subgroup of $H$
  such that $H$ acts on $Y$, that $G$ preserves $V$
  and that moreover
  $G=\{g\in H : g(V)=V\}$.

  Then the morphism $[V/G]\to [Y/H]$ of
  geometric quotients is bimeromorphic
  onto its image.
  \begin{proof} We can replace $Y$ by the largest open subspace
    $Y^0$ of $Y$ such that $Y^0$ is preserved
    by $H$ and $Y^0\to [Y^0/H]$ is {\'e}tale. Then the result
    is clear.
  \end{proof}
\end{lemma}

\begin{lemma} A generic pair $(C,Z)$
  has no automorphisms.
  \begin{proof} We consider the various values
    of $q$ separately.
    \begin{enumerate}
      \item $q\ge 3$. A generic curve
        has no automorphisms.
      \item $q=2$. Then $Z$ is
        invariant under the hyperelliptic involution
        of $C$.
        In this case the dimension of the set $\Sigma$
        of pairs
        $(C,Z)$
    is at most $\dim\sM_2+N/2$, which
    is less than $\dim\sje$.
  \item $q=1$. Either $Z$ is invariant
    under an involution, and then we argue
    as in the case where $q=2$, or $Z$ is invariant
    under a translation. In this case $Z$ is determined
    by any one of its points and
    $\dim\Sigma\le 2$.
  \item $q=0$. Suppose $1\ne g\in PGL_2(\C)$.
    If $g$ fixes one point then $\dim\Sigma\le 2$
    and if $g$ fixes $2$ points then $\dim\Sigma\le 3$.
  \end{enumerate}
    \end{proof}
\end{lemma}

We now take $V$ to be a miniversal deformation
space of a surface $\tX$ with $r=r(q)$ fibres of type $I_2$ and $Y$
a germ of the domain $D$ at the period point of $\tX$.

Suppose that $H\subset O(H^2_{prim}(\tX,\Z))$
is the automorphism group of the polarized
Hodge structure of the $\tX$ and that $G\subset H$
is the subgroup consisting of those elements of $H$ that
preserve the image of $V$ in $Y$. Since $H$ is finite,
we can identify the germs $Y$ and $V$, which are smooth,
with their tangent spaces. In turn,
$V=H^{2,0}(\tX)^\vee\otimes H^{1,1}_{prim}(\tX)$
and $H$ acts on both components of this tensor product.
Let $W=W(A_1^r)$ denote the Weyl group generated
by reflexions $\s_{[\eta_a]}$ in the classes $[\eta_a]$
where $a$ lies over $j=\infty$;
then $\s_{[\eta_a]}$ lies in $H$,
so that $W\subset H$. Also $(\pm 1)$ is a subgroup
of $H$ and acts trivially on $D$.

\begin{proposition}\label{pm} $G= (\pm 1)\times W$.
  \begin{proof} Let $g\in G$.
    We know that the subspace $V$ of $Y$
    is based by the rank $1$ tensors
    $\omega_a^\vee\otimes[\eta_a]$ and that
    these are the only rank $1$ tensors in $V$.
    Therefore $g$ permutes the lines
    $\C\omega_a^\vee$ and so
    acts on the pair $(C,Z)$,
    where $C$ is embedded in $\P^{h-1}$
    via $\vert K_C+L\vert$. This action
    is trivial, by the previous lemma,
    and so $g$ acts as a scalar on $H^{2,0}(\tX)$.
    We can then take this scalar to be $1$.
    Then $g$ fixes each line $\C[\eta_a]$.
    The classes $[\eta_a]$ are orthogonal
    and can be normalized by imposing the
    condition that $[\eta_a]^2=-2$ for all $a\in Z$.
    This normalization is unique up to a choice of
    sign for each $a$.
    Since $\sZ_r$ is irreducible, one choice
    of sign, for a point in $Z_r$, determines
    all the other signs, except when
    $a$ lies over $j=\infty$. This last ambiguity
    is exactly taken care of by the Weyl group
    and so $G\subset (\pm 1)\times W$.

    It is clear that $(\pm 1)\times W\subset G$.
  \end{proof}
\end{proposition}

\begin{theorem} (Generic Torelli) The period map
  $[per]:[\sje]\to [D/\G]$ is birational
  onto its image.
  \begin{proof} The theorem follows from
    Corollary \ref{corollary} and the
    fact that, at the level of miniversal
    deformation spaces, the morphism
    $\tsje\to\sje$ is a quotient map
    by the relevant Weyl group. We conclude
    via Proposition \ref{pm} and the fact that
    $(\pm 1)$ acts trivially.
  \end{proof}
\end{theorem}
\end{section}
\begin{section}{The variational Schottky problem for Jacobian surfaces}\label{schottky}
  The Schottky problem is that of determining the image of a
  moduli space under a period map. As explained by Donagi \cite{D},
  there is a variational approach to this. For curves of genus $q$
  that are neither hyperelliptic, trigonal nor plane quintics,
  his approach leads (see p. 257 of \cite{D}) to the statement that
  the image of the variational period map lies in the
  Grassmannian that parametrizes $3q-3$-dimensional quotient spaces $W$ of
  $\Symm^2V$, where $V$ is a fixed $q$-dimensional vector space,
  and the kernel of $\Symm^2V\to W$ defines a smooth linearly normal
  curve $C$ of genus $q$ in $\P(V)$. (It follows from this that the
  embedding $C\inj\P(V)$ can be identified with the canonical embedding
  of $C$.)
  
  For Jacobian elliptic surfaces of geometric genus $h$ and
  irregularity $q$ we get something equally concrete.
    
  Let $\sZ\to\sJE^{gen}$ be the universal ramification locus,
  of degree $10+8(1-q)$ over $\sJE^{gen}$.
  The image of $\sZ$
  lies in a tensor product
  $U\otimes V$ where $\dim U=h$ and $\dim V=10h+8(1-q)$.
  Projecting to $\P(U)=\P^{h-1}$ leads
    to the following variational partial solution to the Schottky
  problem.
  The solution is only partial because
  this projection factors through the
  quotient stack $\sJE^{gen}/\GG_m$.

  Recall that $\sJE^{gen}$
  can be described as follows.

  Suppose that $\sP=\sP^{h+1-q}\to\sM_q$
is the universal Picard variety of degree $h+1-q$ line bundles
$L$ on a curve $C$ of genus $q$. Then writing
the equation of $X$ in affine Weierstrass
form, namely as
$$y^2=4x^3-g_4x-g_6,$$
shows that
$g_n\in H^0(C,nL)$ and $\sJE^{gen}$ is birationally equivalent to
a $B(\Z/2)$-gerbe over a bundle over $\sP$ whose fibre is the quotient stack
$$(H^0(C,4L)\oplus H^0(C,6L)-\{(0,0)\})/\GG_m.$$
The action of $\GG_m$ on $\bsEll$ leads to an action of $\GG_m$
on $\sJE^{gen}$
and the quotient $\sJE^{gen}/\GG_m$ is birationally equivalent
to a $B(\Z/2)$-gerbe over
the universal $\vert 4L\vert\times\vert 6L\vert$-bundle
over $\sP$.
  
  \begin{theorem} The image of the variational period
    map for Jacobian elliptic surfaces lies in the locus
    $\sV=\sV_{h,q}$
    of zero-cycles $Z$ in $\P^{h-1}$ such that
    \begin{enumerate}
    \item $\deg Z=10h+8(1-q)$;
    \item the intersection of the quadrics through $Z$
      is a curve $C$ of genus $q$ and degree $h+q-1$;
    \item the divisor $Z$ on $C$ is linearly equivalent
      to $10L+K_C$ and the hyperplane class $H$ on $C$
      is $H= L+K_C$. So $Z\sim 10H-9K_C$.
    \end{enumerate}
    \begin{proof} This follows at once from the results
      of the previous section.
    \end{proof}
  \end{theorem}
  We can make this more precise. First, recall the idea of \emph{transvectants}:
  if $N$ is a line bundle on a variety $V$ over a field $k$ and $m,n\in\Z$,
  then there is a homomorphism
  $$N^{\otimes m}\otimes_kN^{\otimes n}\to\Omega^1_V\otimes_{\sO_V} N^{\otimes m+n}$$
  of sheaves on $V$
  defined, in terms of a local generator $s$ of $N$, by
  $$fs^{\otimes m}\otimes_k gs^{\otimes n}\mapsto
  (m f\, dg-n g\, df)\otimes_{\sO_V} s^{\otimes m+n}.$$

  Suppose that $V$ is
  projective; then at the level of global sections
  this defines a $k$-linear homomorphism
  $$H^0(V,N^{\otimes m})\otimes_k H^0(V,N^{\otimes n})\to
  H^0(V,\Omega^1_V\otimes_{\sO_V} N^{\otimes m+n})$$
  of finite-dimensional vector spaces.
  This morphism has several names, depending on the context; for example,
  the first transvectant, the first Ueberschiebung,
  the Jacobian determinant,
  the Gauss map, the Wahl map, the first Rankin--Cohen bracket.
  Its relevance for us lies in the case where $V=C$, $N=L$, $m=4$ and
  $n=6$.

  \begin{lemma} The ramification divisor
    $Z$, which is a point in $\vert 10L+K_C\vert$,
    lies in the image of the projectivized
    first transvectant, which is a rational bilinear map
    $$V_L:\vert 4L\vert\times\vert 6L\vert -\ -\to \vert 10L+K_C\vert.$$
    \begin{proof} Consider the elliptic surface given, in affine terms,
      by the equation $y^2=4x^3-g_4x-g_6,$ where $g_n$ is a
      section of $L^{\otimes n}$. Then the $j$-invariant
      is a fractional linear function
      of the quantity $j_1=g_4^3/g_6^2$, and the lemma follows from
      calculating the zero locus of
      $d j_1/dz$ when $z$ is a local co-ordinate on $C$.
    \end{proof}
  \end{lemma}

  \begin{corollary}\label{6.3} $\sZ$ is irreducible.
    \begin{proof} $\sZ$ is dominated by
      $\vert 6L\vert\times\vert 4L\vert$.
    \end{proof}
  \end{corollary}
  
  Let $V(L)$ denote the image of the projectivized
  first transvectant.

  \begin{theorem} The image of the variational period
    map for Jacobian elliptic surfaces equals the locus $\sT$
    of triples $(C,L,Z)$ where
    \begin{enumerate}
    \item $C\in\sM_q$;
    \item $L\in\Pic^{h+1-q}_C;$
    \item $Z\in V(L)$.
    \end{enumerate}
    This locus is irreducible.
    \begin{proof} The only things left to do 
      is to observe that $\sT$ is irreducible.
      But $\sT$ is dominated by $\sZ$.
    \end{proof}
    \end{theorem}
 \end{section}
\begin{section}{Multiple surfaces}
  An elliptic surface $f:X\to C$ is \emph{multiple} if it is
  not simple; that is, if it has multiple fibres. As far as
  the period map is concerned there is not much to be said
  about these. Kodaira proved that, given $f:X\to C$,
  there is a simple elliptic surface $g:Y\to C$ such that
  $X$ is obtained from $Y$ by logarithmic transformations; in
  particular, there is a
  finite subset $S$ of $C$ such that, if $C_0=C-S$,
  $X_0=f^{-1}(C_0)$ and $Y_0=g^{-1}(C_0)$, then
  $X_0$ and $Y_0$ are isomorphic relative to $C_0$,
  but are not usually bimeromorphic.
  In this situation the sheaves
  $f_*\omega_{X/C}$ and $g_*\omega_{Y/C}$ are isomorphic
  (\cite{Sch}, p. 234).
  It follows that $f_*\Omega^2_X$
  and $g_*\Omega^2_Y$ are isomorphic
  and then that the weight 2 Hodge structures on $X$
  and $Y$ are isomorphic.

  We deduce that
  the period map can detect neither the
  presence nor the location of multiple fibres on an elliptic surface.
\end{section}

\end{document}